\documentclass[12pt]{article}
\usepackage{latexsym}
\usepackage{amsmath}
\usepackage{amssymb}
\usepackage{amsfonts}
\usepackage{graphics,graphicx}
\usepackage{mathrsfs}
\usepackage{cite}
\usepackage{color}
\usepackage[colorlinks=true,linkcolor=blue]{hyperref}
\usepackage{tikz}

\numberwithin{equation}{section}

\newtheorem{theorem}{Theorem}[section]
\newtheorem{lemma}{Lemma}[section]
\newtheorem{pr}{Proposition}[section]

\newtheorem{definition}{Definition}[section]
\newtheorem{remark}{Remark}[section]

\newcommand{\eproof}{{\mbox{\ }~\hfill
\mbox{\large $\Box$} \par \vskip 10pt}}
\newcommand{\pf}{\noindent{\bf Proof}}
\newcommand{\eps}{\varepsilon}
\renewcommand{\l}{\ell}

\newcommand{\supp}{\mbox{\rm supp}}
\newcommand{\dom}{\Omega}
\newcommand{\sig}{\Sigma}

\newcommand{\el}{\mathcal{L}}

\newcommand{\R}{\mathbb R}

\title{Propagation of smallness and size estimate in the second order elliptic equation with discontinuous complex Lipschitz conductivity}
\author{Elisa~Francini\thanks{Universit\`a di Firenze, Italy. Email: elisa.francini@unifi.it}
\qquad Sergio~Vessella\thanks{Universit\`a di Firenze, Italy. Email: sergio.vessella@unifi.it}
\qquad Jenn-Nan~Wang\thanks{National Taiwan University, Taiwan. Email: jnwang@math.ntu.edu.tw}}

\date{\today}

\begin{document}
\maketitle

\begin{abstract}
In this paper, we would like to derive three-ball inequalities and propagation of smallness for the complex second order elliptic equation with discontinuous Lipschitz coefficients. As an application of such estimates, we study the size estimate problem by one pair of Cauchy data on the boundary. The main ingredient in the derivation of three-ball inequalities and propagation of smallness is a local Carleman proved in our recent paper \cite{FVW20}.

\end{abstract}

\section{Introduction}

In our recent paper \cite{FVW20}, we derived a Carleman estimate for the second order elliptic equation with piecewise complex-valued Lipschitz coefficients. The theme of this paper is to prove some interesting results based on the Carleman estimate obtained in \cite{FVW20}. The ultimate goal is to give upper and lower bounds of the size of the inclusion embedded inside of a conductive body with discontinuous complex conductivity by only one pair of boundary measurements. A typical application of this study is to estimate the size of a cancerous tumor inside an organ by the electric impedance tomography (EIT).  

For the conductivity equation with piecewise real Lipschitz coefficients, the same size estimate problem was considered in \cite{flvw}. We want to point out that, in many real world problems, the case of complex-valued coefficients arises naturally. The modeling of the current flows in biological tissues or the propagation of the electromagnetic waves in conductive media are typical examples. In these cases, the conductivities are complex-valued functions. On the other hand, in some situations, the conductivities are not continuous functions. In the human body, different organs have different conductivities. For instance, the conductivities of heart, liver, intestines are 0.70
(S/m), 0.10 (S/m), 0.03 (S/m), respectively. Therefore, to model the current flow in the human body, it is more reasonable to consider an anisotropic complex-valued conductivity with jump-type discontinuities \cite{mph06}.      

In the size estimate problem studied in \cite{flvw}, the essential tool is a three-region inequality which is obtained by applying the Carleman estimate for the second order elliptic equation with piecewise real Lipschitz coefficients derived in \cite{dflvw}. Since we have the similar Carleman estimate available for the case of piecewise complex-valued Lipschitz coefficients, we can proceed the method used in \cite{flvw} to prove the three-region inequality. In treating the size estimate problem, the three-region inequality is enough since one only needs to cross the interface once in propagating the information in the interior to the boundary. However, the three-region inequality is inconvenient in deriving the general propagation of smallness. Therefore, in this paper, we want to derive the usual three-ball inequality for the complex second order elliptic operator even when the coefficients are piecewise Lipschitz. We will follow the ideas outlined in \cite{CW} where the three-ball inequality was proved for the real second order elliptic operator with piecewise Lipschitz coefficients. We then apply the three-ball inequality to derive a general propagation of smallness for the second order elliptic equation with piecewise complex Lipschitz coefficients.

We would like to raise an issue in the investigation of the size estimate problem when the background medium is complex valued. Following the method used in \cite{ars}, an important step is to derive certain energy inequalities controlling the power gap. It was noted in \cite{bfv} that when the current flows inside and outside of the inclusion obey the usual Ohm's law (the relation between current and voltage is linear) and the imaginary part of the conductivity outside of the inclusion is a nonzero variable function, energy inequalities \eqref{eineq1}, \eqref{eineq2} are not likely to hold.  Precisely, in this case, an example with $\delta W=0$, but $|D|\ne 0$, is constructed in \cite{bfv}. On the other hand, if both conductivities inside and outside of the inclusion are complex constants or the conductivity outside of the inclusion is real valued, then energy inequalities \eqref{eineq1}, \eqref{eineq2} were obtained in \cite{bfv}.  

A key observation found in \cite{CNW} is that if the current flow inside of the inclusion obeys certain nonlinear Ohm's law, we can restore the energy inequalities \eqref{eineq1}, \eqref{eineq2}. In particular, our size estimate result applies to the case of a non-chiral medium with a chiral inclusion having real valued chirality.

The paper is organized as follows. In Section~\ref{pre}, we introduce some notations and state the Carleman estimate proved in \cite{FVW20}. In Section~\ref{3regions}, we plan to prove a three-region inequality across the interface based on the Carleman estimate given in Section~\ref{pre}. We then combine the classical three-ball inequality and the three-region inequality to derive a three-ball inequality in Section~\ref{smallness}. There, we also prove the propagation of smallness. Finally, we study the size estimate problem in Section~\ref{size-est}.

\section{Notations and Carleman estimate}\label{pre}

In this section, we will state the Carleman {{estimate}} proved in \cite{FVW20} where the interface is assumed to be flat. Since our Carleman estimate is local near any point at the interface, for {a} general $C^{1,1}$ interface, it can be flatten by a suitable change of coordinates. Moreover, the transformed coefficients away from the interface remain Lipschitz. Define $H_{\pm}=\chi_{\mathbb{R}^n_{\pm}}$ where $\mathbb{R}^n_{\pm}=\{(x',x_n)\in \mathbb{R}^{n-1}\times\mathbb{R}|x_n\gtrless0\}$ and $\chi_{\mathbb{R}^n_{\pm}}$ is the characteristic function of $\mathbb{R}^n_{\pm}$.
In places we will use equivalently the symbols $\partial$, $\nabla$ and $D=-i\nabla$
to denote the gradient of a function and we will add the index $x'$ or $x_n$
to denote gradient in $\mathbb R^{n-1}$ and the derivative with respect to $x_n$
respectively. We further denote $\partial_\ell=\partial/\partial x_\ell$, $D_\ell=-i\partial_\ell$, and $\partial_{\xi_\ell}=\partial/\partial {\xi_\ell}$. 

Let $u_\pm\in C^\infty(\R^n)$. We define
\begin{equation*}\label{1.030}
u=H_+u_++H_-u_-=\sum_\pm H_{\pm}u_{\pm},
\end{equation*}
hereafter, we denote $\sum_\pm a_\pm=a_++a_-$,  and 
\begin{equation*}
\mathcal{L}(x,D)u:=\sum_{\pm}H_{\pm}{\rm div}(A_{\pm}(x)\nabla u_{\pm}),
\end{equation*}
where
\begin{equation}\label{7.2}
A_{\pm}(x)=\{a^{\pm}_{\l j}(x)\}^n_{\l ,j=1}=\{a^{\pm}_{\l j}(x',x_n)\}^n_{\l ,j=1},\quad x'\in \mathbb{R}^{n-1},x_n\in \mathbb{R}
\end{equation}
is a Lipschitz symmetric matrix-valued function. Assume that
\begin{equation}\label{symm0}
a_{\l j}^{\pm}(x)=a_{j\l}^{\pm}(x),\quad\forall\;\;\l, j=1,\cdots,n,
\end{equation}
and furthermore
\begin{equation}\label{complex0}
a_{\l j}^{\pm}(x)=M_{\l j}^{\pm}(x)+i \gamma N_{\l j}^{\pm}(x),
\end{equation}
where $(M_{\l j}^{\pm})$ and $(N_{\l j}^{\pm})$ are real-valued matrices and $\gamma> 0$. We further assume that there exists $\lambda_0>0$ such that for all $\xi\in\R^n$ and $x\in \R^n$ we have
\begin{equation}\label{elliptic10}
\lambda_0|\xi|^2\le M^{\pm}(x)\xi\cdot\xi\le\lambda_0^{-1}|\xi|^2
\end{equation}
and
\begin{equation}\label{elliptic20}
\lambda_0|\xi|^2\le N^{\pm}(x)\xi\cdot\xi\le\lambda_0^{-1}|\xi|^2.
\end{equation}
In the paper, we consider Lipschitz coefficients $A_{\pm}$, i.e., there exists a constant $M_0>0$ such that
\begin{equation}\label{7.4}
|A_{\pm}(x)-A_{\pm}(y)|\leq M_0|x-y|.
\end{equation}
To treat the transmission conditions, we write
\begin{equation}\label{7.5}
h_0(x'):=u_+(x',0)-u_-(x',0),\ \forall\,\; x'\in \mathbb{R}^{n-1},
\end{equation}
\begin{equation}\label{7.6}
h_1(x'):=A_+(x',0)\nabla u_+(x',0)\cdot \nu-A_-(x',0)\nabla u_-(x',0)\cdot \nu,\ \forall\,\; x'\in \mathbb{R}^{n-1},
\end{equation}
where $\nu=e_n$.

Let us now introduce  the weight function. Let $\varphi$ be
\begin{equation}\label{2.1}
\varphi(x_n)=
\begin{cases}
\begin{array}{l}
\varphi_+(x_n):=\alpha_+x_n+\beta x_n^2/2,\quad x_n\geq 0,\\
\varphi_-(x_n):=\alpha_-x_n+\beta x_n^2/2,\quad x_n< 0,
\end{array}
\end{cases}
\end{equation}
where $\alpha_+$, $\alpha_-$ and $\beta$ are positive numbers which will be determined later. In what follows we denote by $\varphi_{+}$ and $\varphi_{-}$ the restriction of the weight function $\varphi$ to $[0,+\infty)$ and to $(-\infty,0)$ respectively. We use similar notation for any other weight functions. For any $\varepsilon>0$ let 
\begin{equation}\label{psi}
\psi_{\varepsilon}(x):=\varphi(x_n)-\frac{\varepsilon}{2}|x'|^2,
\end{equation}
and let 
\begin{equation}\label{wei}
\phi_{\delta}(x):=\psi_{\delta}(\delta^{-1}x),\quad\delta>0.
\end{equation}

For a function $h\in L^2(\mathbb{R}^{n})$, we define
\begin{equation*}
\hat{h}(\xi',x_n)=\int_{\mathbb{R}^{n-1}}h(x',x_n)e^{-ix'\cdot\xi}\,dx',\quad \xi'\in \mathbb{R}^{n-1}.
\end{equation*}
As usual we denote by $H^{1/2}(\mathbb{R}^{n-1})$ the space of the functions $f\in L^2(\mathbb{R}^{n-1})$ satisfying
$$\int_{\mathbb{R}^{n-1}}|\xi'||\hat{f}(\xi')|^2d\xi'<\infty,$$
with the norm
\begin{equation}\label{semR}
\|f\|^2_{H^{1/2}(\mathbb{R}^{n-1})}=\int_{\mathbb{R}^{n-1}}(1+|\xi'|^2)^{1/2}|\hat{f}(\xi')|^2d\xi'.
\end{equation}
Moreover we define
$$[f]_{1/2,\mathbb{R}^{n-1}}=\left[\int_{\mathbb{R}^{n-1}}\int_{\mathbb{R}^{n-1}}\frac{|f(x)-f(y)|^2}{|x-y|^n}dydx\right]^{1/2},$$
and recall that there is a positive constant $C$, depending only on $n$, such that
\begin{equation*}
C^{-1}\int_{\mathbb{R}^{n-1}}|\xi'||\hat{f}(\xi')|^2d\xi'\leq[f]^2_{1/2,\mathbb{R}^{n-1}}\leq C\int_{\mathbb{R}^{n-1}}|\xi'||\hat{f}(\xi')|^2d\xi',
\end{equation*}
so that the norm \eqref{semR} is equivalent to the norm $\|f\|_{L^2(\mathbb{R}^{n-1})}+[f]_{1/2,\mathbb{R}^{n-1}}$. We use the letters $C, C_0, C_1, \cdots$ to denote constants. The value of the constants may change from line to line, but it is always greater than $1$.

We will denote by $B'_r(x')$ the $(n-1)$-ball centered at $x'\in \mathbb{R}^{n-1}$ with radius $r>0$. Whenever $x'=0$ we denote $B'_r=B'_r(0)$. Likewise, we denote $B_r(x)$ be the $n$-ball centered at $x\in\R^n$ with radius $r>0$ and $B_r=B_r(0)$.
\bigskip

\begin{theorem}\label{thm8.2}
Let $A_{\pm}(x)$ satisfy \eqref{7.2}-\eqref{7.4}. There exist $\alpha_+,\alpha_-,\beta, \delta_0, r_0, \gamma_0$, $\tau_0$, $C$ depending on $\lambda_0, M_0$ such that if $\gamma\le\gamma_0$, $\delta\le\delta_0$ and $\tau\geq \tau_0$, then
\begin{equation}\label{8.24}
\begin{aligned}
&\sum_{\pm}\sum_{k=0}^2\tau^{3-2k}\int_{\mathbb{R}^n_{\pm}}|D^k{u}_{\pm}|^2e^{2\tau\phi_{\delta,\pm}(x',x_n)}dx'dx_n+\sum_{\pm}\sum_{k=0}^1\tau^{3-2k}\int_{\mathbb{R}^{n-1}}|D^k{u}_{\pm}(x',0)|^2e^{2\phi_\delta(x',0)}dx'\\
&+\sum_{\pm}\tau^2[e^{\tau\phi_{\delta}(\cdot,0)}u_{\pm}(\cdot,0)]^2_{1/2,\mathbb{R}^{n-1}}+\sum_{\pm}[D(e^{\tau\phi_{\delta,\pm}}u_{\pm})(\cdot,0)]^2_{1/2,\mathbb{R}^{n-1}}\\
\leq &C\left(\sum_{\pm}\int_{\mathbb{R}^n_{\pm}}|\mathcal{L}(x,D)(u_{\pm})|^2\,e^{2\tau\phi_{\delta,\pm}(x',x_n)}dx'dx_n+[e^{\tau\phi_\delta(\cdot,0)}h_1]^2_{1/2,\mathbb{R}^{n-1}}\right.\\
&\left.+[D_{x'}(e^{\tau\phi_\delta}h_0)(\cdot,0)]^2_{1/2,\mathbb{R}^{n-1}}+\tau^{3}\int_{\mathbb{R}^{n-1}}|h_0|^2e^{2\tau\phi_\delta(x',0)}dx'+\tau\int_{\mathbb{R}^{n-1}}|h_1|^2e^{2\tau\phi_\delta(x',0)}dx'\right).
\end{aligned}
\end{equation}
where $u=H_+u_++H_-u_-$,  $u_{\pm}\in C^\infty(\mathbb{R}^{n})$ and ${\rm supp}\, u\subset B'_{\delta r_0}\times[-\delta r_0,\delta r_0]$, and $\phi_\delta$ is given by \eqref{wei}.
\end{theorem}
\begin{remark}
In view of the proof of Theorem~\ref{thm8.2} in {\rm\cite{FVW20}}, the coefficients $\alpha_+, \alpha_-$ are required to satisfy
\begin{equation}\label{k0}
\frac{\alpha_+}{\alpha_-}\ge \kappa_0>1,
\end{equation}
where $\kappa_0$ is general constant depending on the values of $A_{\pm}(0)$ at the interface.
\end{remark}
\begin{remark}
It is clear that \eqref{8.24} remains valid if can add lower order terms\\ $\sum_{\pm}H_{\pm}\left(W\cdot\nabla u_{\pm}+Vu_{\pm}\right)$, where $W,V$ are bounded functions, to the operator ${\mathcal L}$. That is, one can substitute
\begin{equation}\label{8.242}
\mathcal{L}(x,D)u=\sum_{\pm}H_{\pm}{\rm div}(A_{\pm}(x)\nabla u_{\pm})+\sum_{\pm}H_{\pm}\left(W\cdot\nabla u_{\pm}+Vu_{\pm}\right)
\end{equation}
in \eqref{8.24}.
\end{remark}

\section{Three-region inequalities}\label{3regions}

Based on the Carleman estimate given in Theorem~\ref{thm8.2}, we will derive three-region inequalities across the interface $x_n=0$. Here we consider $u=H_+u_++H_-u_-$ satisfying
\begin{equation}\label{eq101}
\mathcal{L}(x,D)u=0\quad\mbox{in}\quad\R^n,
\end{equation}
where ${\mathcal L}$ is given in \eqref{8.242} and
\[
\|W\|_{L^\infty(\R^n)}+\|V\|_{L^\infty(\R^n)}\le\lambda_0^{-1}.
\]
Now all coefficients $\alpha_\pm$, $\beta$, $\delta_0$, $\gamma_0$, $r_0$, $\tau_0$ have been determined in Theorem~\ref{thm8.2}.

\begin{theorem}\label{thm9.1}
Let $u$ be a solution of \eqref{eq101} and $A_{\pm}(x)$ satisfy \eqref{7.2}-\eqref{7.4} with $h_0=h_1=0$. Moreover, the constant $\gamma$ in \eqref{complex0} satisfies $\gamma\le\gamma_0$ with $\gamma_0$ given in Theorem~\ref{thm8.2}. Then there exist $C$ and ${R}$, depending only on $\lambda_0, M_0, n$, such that if $\ 0<R_1,R_2\leq R$, then
\begin{equation}\label{9.1}
\int_{U_2}|u|^2 dx \le (e^{\tau_0R_2}+CR_1^{-4}) \left(\int_{U_1}|u|^2dx\right)^{\frac{R_2}{2R_1+3R_2}}\left(\int_{U_3}|u|^2dx\right)^{\frac{2R_1+2R_2}{2R_1+3R_2}},
\end{equation}
where 
\[
\begin{aligned}
&U_1=\left\{z\geq-4R_2,\,\frac{R_1}{8a}<x_n<\frac{R_1}{a}\right\},\\
&U_2=\left\{-R_2\leq z\leq\frac{R_1}{2a},\,x_n<\frac{R_1}{8a}\right\},\\
&U_3=\left\{z\geq-4R_2,\,x_n<\frac{R_1}{a}\right\},
\end{aligned}
\]
$a=\alpha_+/\delta$,
\begin{equation}\label{zxy}
z(x)=\frac{\alpha_-x_n}{\delta}+\frac{\beta x_n^2}{2\delta^2}-\frac{|x'|^2}{2\delta},
\end{equation}
and any $\delta\le\delta_0$.
\end{theorem}
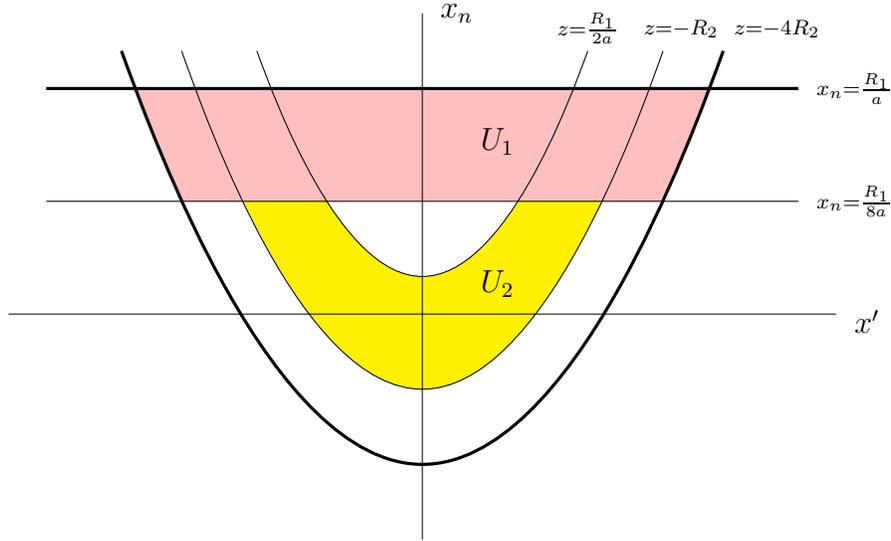
\begin{figure}[ht]
\centering
\begin{tikzpicture}
\begin{scope}
\clip (-5,3) rectangle (5,1.5);
\fill[pink] (-4,3.5) parabola bend (0,-2) (4,3.5);
\end{scope}
\draw[line width=1.2pt, black] (-4,3.5) parabola bend (0,-2) (4,3.5);
\begin{scope}
\clip (-5,1.5) rectangle (5,-1.2);
\fill[yellow] (-3.2,3.5) parabola bend (0,-1) (3.2,3.5);
\fill[white] (-2.2,3.5) parabola bend (0,0.5) (2.2,3.5);
\end{scope}
\draw[line width=0.2pt]  (5.5,0) -- (-5.5,0);
\node[right] at (0.1,4) {$x_n$}; 
\draw[line width=0.2pt]  (0,4) -- (0,-3);
\node[right] at (5.6,-0.1) {$x'$};
\node[right] at (4.0,3.8) {$\scriptstyle z=-4R_2$}; 
\draw (-3.2,3.5) parabola bend (0,-1) (3.2,3.5);
\node[right] at (2.8,3.8) {$\scriptstyle z=-R_2$};
\draw (-2.2,3.5) parabola bend (0,0.5) (2.2,3.5);
\draw[line width=1.2pt, black] (-5,3) -- (5,3);
\draw (-5,1.5) -- (5,1.5);
\node at (2.2,3.8) {$\scriptstyle z=\frac{R_1}{2a}$};
\node[right] at (5.1,3) {$\scriptstyle x_n=\frac{R_1}{a}$};
\node[right] at (5.1,1.5) {$\scriptstyle x_n=\frac{R_1}{8a}$};
\node at (1,2.3) {$U_1$};
\node at (1,0.4) {$U_2$};
\end{tikzpicture}
\caption{$U_1$ and $U_2$ are shown in pink and yellow, respectively. $U_3$ is the region enclosed by black boundaries. Note that since $z$ is hyperbolic, there are parts similar to $U_2$ and $U_3$ lying below $x_n<-\alpha_-\delta/\beta$. Here we are only interested in the solution near $x_n=0$. Thus we consider the cut-off function relative to $U_2$ and $U_3$ as in the figure.}
\label{fig1}
\end{figure}

\pf. Here we adopt the proof given in \cite{flvw}. To apply the estimate \eqref{8.24}, we need to ensure that $u$ satisfies the support condition. Let $r>0$ be chosen satisfying
\begin{equation}\label{1.5}
r\le\min\left\{r_0^2,\frac{13\alpha_-}{8\beta},\frac{2\delta r_0}{19\alpha_-+8\beta}\right\}.
\end{equation}
We then set
\[
R=\frac{\alpha_- r}{16}.
\]
It follows from \eqref{1.5} that
\begin{equation}\label{bee}
R\le\frac{13\alpha_-^2}{128\beta}.
\end{equation}
Given $0<R_1<R_2\le R$. Let $\vartheta_1(t)\in C^{\infty}_0 ({\mathbb R})$ satisfy $0\le\vartheta_1(t)\leq 1$ and
\begin{equation*}
\vartheta_1 (t)=\left\{
\begin{aligned}
&1,\quad t>-2R_2,\\
&0,\quad t\leq -3R_2.
\end{aligned}\right.
\end{equation*}
Also, define $\vartheta_2(t)\in C^{\infty}_0 ({\mathbb R})$ satisfying $0\le\vartheta_2(t)\leq 1$ and
\begin{equation*}
\vartheta_2 (t)=\left\{
\begin{aligned}
&0,\quad t\geq \frac{R_1}{2a},\\
&1,\quad t<\frac{R_1}{4a}.
\end{aligned}\right.
\end{equation*}
Finally, we define $\vartheta(x)=\vartheta(x',x_n)=\vartheta_1(z(x))\vartheta_2(x_n)$, where $z$ is defined by \eqref{zxy}. 

We now check the support condition for $\vartheta$. From its definition, we can see that $\supp\,\vartheta$ is contained in
\begin{equation}\label{2.2}
\left\{
\begin{aligned}
&z(x)=\frac{\alpha_-x_n}{\delta}+\frac{\beta x_n^2}{2\delta^2}-\frac{|x'|^2}{2\delta}>-3R_2,\\
& x_n<\frac{R_1}{2a}.
\end{aligned}\right.
\end{equation}
In view of the relation
\[
\alpha_+>\alpha_-\;\;(\mbox{see}\;\eqref{k0})\quad\mbox{and}\quad a=\frac{\alpha_+}{\delta},
\]
we have that
\begin{equation*}
\frac{R_1}{2a}<\frac{\delta}{2\alpha_-}\cdot R_1<\frac{\delta}{\alpha_-}\cdot\frac{\alpha_-r}{16}<\delta r,
\end{equation*}
i.e., $x_n<\delta r\le\delta r^2_0\le\delta r_0$. Next, we  observe that
\[
-3R_2>-3R=-\frac{3\alpha_-r}{16}>\frac{\alpha_-}{\delta}(-\delta r)+\frac{\beta}{2\delta^2}(-\delta r)^2,
\]
which gives $-\delta r_0<-\delta r<x_n$ due to  \eqref{1.5}. Consequently, we verify that $|x_n|<\delta r<\delta r_0$. One the other hand, from the first condition of \eqref{2.2} and \eqref{1.5}, we see that
\begin{equation*}
\begin{aligned}
\frac{|x'|^2}{2\delta}&<3R_2+\frac{\alpha_-x_n}{\delta}+\frac{\beta x_n^2}{2\delta^2}\leq\frac{3\alpha_-r}{16}+\frac{\alpha_-}{\delta}\cdot \delta r+\frac{\beta}{2\delta^2}\cdot \delta^2r^2\\
&\le\frac{19\alpha_-+8\beta}{16}r\le\frac{19\alpha_-+8\beta}{16}r_0^2\leq\frac{\delta}{8}r_0^2,
\end{aligned}
\end{equation*}
which gives $|x'|<\delta r_0/2$.

Since $h_0=0$, we have that
\begin{equation}\label{8.90}
\vartheta(x',0)u_+(x',0)-\vartheta(x',0)u_-(x',0)=0,\;\forall\; x'\in\R^{n-1}.
\end{equation}
Applying \eqref{8.24} to $\vartheta u$ and using \eqref{8.90} yields
\begin{equation}\label{9.2}
\begin{aligned}
&\sum_{\pm}\sum_{|k|=0}^2\tau^{3-2|k|}\int_{\mathbb{R}^n_{\pm}}|D^k(\vartheta u_{\pm})|^2e^{2\tau\phi_{\delta,\pm}(x',x_n)}dx'dx_n\\
\le &C\sum_{\pm}\int_{\mathbb{R}^n_{\pm}}|{\mathcal L}(x,D)(\vartheta u_{\pm})|^2\,e^{2\tau\phi_{\delta,\pm}(x',x_n)}dx'dx_n\\
&+C\tau\int_{\mathbb{R}^{n-1}}|A_+(x',0)\nabla (\vartheta u_+(x',0))\cdot \nu-A_-(x',0)\nabla (\vartheta u_-)(x',0)\cdot \nu|^2e^{2\tau\phi_\delta(x',0)}dx'\\
&+C[e^{\tau\phi_\delta(x',0)}\big(A_+(x',0)\nabla (\vartheta u_+)(x',0)\cdot \nu-A_-(x',0)\nabla (\vartheta u_-)(x',0)\cdot \nu\big)]^2_{1/2,\mathbb{R}^{n-1}}.
\end{aligned}
\end{equation}
We now observe that $\nabla \vartheta_1(z)=\vartheta_1'(z)\nabla z=\vartheta_1'(z)(-\frac{x'}{\delta},\frac{\alpha_-}{\delta}+\frac{\beta x_n}{\delta^2})$ and it is nonzero only when 
\[
-3R_2<z<-2R_2.
\]
Therefore, when $x_n=0$, we have
\[
2R_2<\frac{|x'|^2}{2\delta}<3R_2.
\]
Thus, we can see that
\begin{equation}\label{1001}
|\nabla\vartheta(x',0)|^2\le CR_2^{-2}\left(\frac{6R_2}{\delta}+\frac{\alpha_-^2}{\delta^2}\right)\le C R_2^{-2}.
\end{equation}
By $h_0(x')=h_1(x')=0$, \eqref{1001}, and the easy estimate of \cite[Proposition~4.2]{dflvw}, we can estimate
\begin{equation}\label{9.3}
\begin{aligned}
&\tau\int_{\mathbb{R}^{n-1}}|A_+(x',0)\nabla (\vartheta u_+(x',0))\cdot \nu-A_-(x',0)\nabla(\vartheta u_-)(x',0)\cdot \nu|^2e^{2\tau\phi_\delta(x',0)}dx'\\
&+[e^{\tau\phi_\delta(x',0)}\big(A_+(x',0)\nabla(\vartheta u_+)(x',0)\cdot \nu-A_-(x',0)\nabla(\vartheta u_-)(x',0)\cdot \nu\big)]^2_{1/2,\mathbb{R}^{n-1}}\\
\le &\,C R_2^{-2}e^{-4\tau R_2}\left(\tau\int_{\{\sqrt{4\delta R_2}\leq|x'|\leq\sqrt{6\delta R_2}\}}|u_+(x',0)|^2dx'+[u_+(x',0)]^2_{1/2,\{\sqrt{4\delta R_2}\leq|x'|\leq\sqrt{6\delta R_2}\}}\right)\\
&+C\tau^2R_2^{-3}e^{-4\tau R_2}\int_{\{\sqrt{4\delta R_2}\leq|x'|\leq\sqrt{6\delta R_2}\}}|u_+(x',0)|^2dx'\\
\le &\,C\tau^2R_2^{-3}e^{-4\tau R_2}E,
\end{aligned}
\end{equation}
where 
\[
E=\int_{\{\sqrt{4\delta R_2}\leq|x'|\leq\sqrt{6\delta R_2}\}}|u_+(x',0)|^2dx'+[u_+(x',0)]^2_{1/2,\{\sqrt{4\delta R_2}\leq|x'|\leq\sqrt{6\delta R_2}\}}.
\]

Writing out ${\cal L}(x,D)(\vartheta u_{\pm})$ and considering the set where $\nabla\vartheta\neq0$, it is not hard to estimate
\begin{equation}\label{9.4}
\begin{aligned}
&\sum_{\pm}\sum_{|k|=0}^1\tau^{3-2|k|}\int_{\{-2R_2\leq z\leq \frac{R_1}{2a},\,x_n<\frac{R_1}{4a}\}}|D^ku_{\pm}|^2e^{2\tau\phi_{\delta,\pm}(x',x_n)}dx'dx_n\\
\le &\,C\sum_{\pm}\sum_{|k|=0}^1R_2^{2(|k|-2)}\int_{\{-3R_2\leq z\leq-2R_2,\,x_n<\frac{R_1}{2a}\}}|D^ku_{\pm}|^2e^{2\tau\phi_{\delta,\pm}(x',x_n)}dx'dx_n\\
&\,+C\sum_{|k|=0}^1R_1^{2(|k|-2)}\int_{\{-3R_2\leq z,\,\frac{R_1}{4a}< x_n<\frac{R_1}{2a}\}}|D^ku_{+}|^2e^{2\tau\phi_{\delta,+}(x',x_n)}dx'dx_n\\
&\,+C\tau^2R_2^{-3}e^{-4\tau R_2}E\\
\le &\, C\sum_{\pm}\sum_{|k|=0}^1R_2^{2(|k|-2)}e^{-4\tau R_2}e^{2\tau\frac{(\alpha_+-\alpha_-)}{\delta}\frac{R_1}{4a}}\int_{\{-3R_2\leq z\leq-2R_2,\,x_n<\frac{R_1}{4a}\}}|D^{k}u_\pm|^2dx'dx_n\\
&+\sum_{|k|=0}^1R_1^{2(|k|-2)}e^{2\tau\frac{\alpha_+}{\delta}\frac{R_1}{2a}}e^{2\tau\frac{\beta}{2\delta^2}(\frac{R_1}{2a})^2}\int_{\{z\geq-3R_2,\,\frac{R_1}{4a}<x_n<\frac{R_1}{2a}\}}|D^{k}u_+|^2dx'dx_n\\
&\,+C\tau^2R_2^{-3}e^{-4\tau R_2}E.
\end{aligned}
\end{equation}
Let us {{recall}} $U_1=\{z\geq-4R_2,\,\frac{R_1}{8a}<x_n<\frac{R_1}{a}\}$,
$U_2=\{-R_2\leq z\leq\frac{R_1}{2a},\,x_n<\frac{R_1}{8a}\}$. From \eqref{9.4} and interior estimates (Caccioppoli's type inequality), we can derive that

\begin{equation}\label{9.5}
\begin{aligned}
&\tau^3e^{-2\tau R_2}\int_{U_2}|u|^2dx'dx_n{{=\tau^{3}e^{-2\tau R_2}\int_{\{-R_2\leq z\leq\frac{R_1}{2a},\,x_n<\frac{R_1}{8a}\}}|u|^2dx'dx_n}}\\
\leq&\,\sum_{\pm}\tau^{3}\int_{\{-2R_2\leq z\leq\frac{R_1}{2a},\,x_n<\frac{R_1}{4a}\}}| u_{\pm}|^2e^{2\tau\phi_{\delta,\pm}(x',x_n)}dx'dx_n\\
\le &\, C\sum_{\pm}\sum_{|k|=0}^1R_2^{2(|k|-2)}e^{-4\tau R_2}e^{2\tau\frac{(\alpha_+-\alpha_-)}{\delta}\frac{R_1}{4a}}\int_{\{-3R_2\leq z\leq-2R_2,\,x_n<\frac{R_1}{4a}\}}|D^{k}u_\pm|^2dx'dx_n\\
&+\sum_{|k|=0}^1R_1^{2(|k|-2)}e^{2\tau\frac{\alpha_+}{\delta}\frac{R_1}{2a}}e^{2\tau\frac{\beta}{2\delta^2}(\frac{R_1}{2a})^2}\int_{\{z\geq-3R_2,\,\frac{R_1}{4a}<x_n<\frac{R_1}{2a}\}}|D^{k}u_+|^2dx'dx_n\\
&\,+C\tau^2R_2^{-3}e^{-4\tau R_2}E\\
\le &\,CR_1^{-4}e^{-3\tau R_2}\int_{\{-4R_2\leq z\leq-R_2,\,x_n<\frac{R_1}{a}\}}| u|^2dx'dx_n+C\tau^2R_2^{-3}e^{-4\tau R_2}E\\
&\,+CR_1^{-4}e^{(1+\frac{\beta R_1}{4\alpha_-^2})\tau R_1}\int_{\{z\geq-4R_2,\,\frac{R_1}{8a}<x_n<\frac{R_1}{a}\}}| u|^2dx'dx_n\\
\le &CR_1^{-4}\left(e^{2\tau R_1}\int_{U_1}|u|^2dx'dx_n+\tau^2e^{-3\tau R_2}F\right),
\end{aligned}
\end{equation}
where \[F=\int_{{{U_3}}}| u|^2dx'dx_n\] and we used the inequality $\frac{\beta R_1}{4\alpha_-^2}<1$ in view of \eqref{bee}. {{Remark that we estimate $E$ by $F$ using the trace estimate and the interior estimate.}}

Dividing $\tau^3e^{-2\tau R_2}$ on both sides of \eqref{9.5} gives
\begin{equation}\label{9.6}
\int_{U_2}|u|^2{{dx}}\le CR_1^{-4}\left( e^{2\tau (R_1+R_2)}\int_{U_1}|u|^2dx'dx_n+e^{-\tau R_2}F\right).
\end{equation}
We now discuss two cases. If $\int_{U_1} |u|^2dx'dx_n\ne 0$ and
$$e^{2\tau_0 (R_1+R_2)}\int_{U_1}|u|^2dx'dx_n<e^{-\tau_0 R_2}F,$$
then we can choose a $\tau>\tau_0$ so that
$$
e^{2\tau (R_1+R_2)}\int_{U_1}|u|^2dx'dx_n=e^{-\tau R_2}F.
$$
With such $\tau$,  it follows from \eqref{9.6} that
\begin{equation}\label{9.7}
\begin{aligned}
\int_{U_2}|u|^2{{dx}}&\le CR_1^{-4} e^{2\tau (R_1+R_2)}\int_{U_1}|u|^2dx'dx_n\\
&=CR_1^{-4}\left(\int_{U_1}|u|^2dx'dx_n\right)^{\frac{R_2}{2R_1+3R_2}}(F)^{\frac{2R_1+2R_2}{2R_1+3R_2}}.
\end{aligned}
\end{equation}
If $\int_{U_1}|u|^2dx'dx_n= 0$, then letting $\tau\to\infty$ in \eqref{9.6} we have $\int_{U_2}|u|^2dx'dx_n=0$ as well. The
three-regions inequality \eqref{9.1} obviously holds.

On the other hand, if
$$ e^{-\tau_0 R_2}F\leq e^{2\tau_0 (R_1+R_2)}\int_{U_1}|u|^2dx'dx_n,$$
then we have
\begin{equation}\label{9.8}
\begin{aligned}
\int_{U_2}|u|^2 dx'dx_n&\leq \left(F\right)^{\frac{R_2}{2R_1+3R_2}}\left(F\right)^{\frac{2R_1+2R_2}{2R_1+3R_2}}\\
&\leq \exp{(\tau_0R_2)}\left(\int_{U_1}|u|^2dx'dx_n\right)^{\frac{R_2}{2R_1+3R_2}}\left(F\right)^{\frac{2R_1+2R_2}{2R_1+3R_2}}.
\end{aligned}
\end{equation}
Putting together \eqref{9.7}, \eqref{9.8} implies
\begin{equation}\label{9.9}
\int_{U_2}|u|^2 dx'dx_n \le
(\exp{(\tau_0R_2)}+CR_1^{-4})\left(\int_{U_1}|u|^2dx'dx_n\right)^{\frac{R_2}{2R_1+3R_2}}\left(F\right)^{\frac{2R_1+2R_2}{2R_1+3R_2}}.
\end{equation}
\eproof

\section{Propagation of smallness}\label{smallness}

In this section, we will derive a general propagation of smallness for solutions satisfying \eqref{eq101}, ${\cal L}(x,D)u=0$ in $\R^n$, using the ideas given in \cite{CW}. For the region away from the interface, classical three-ball inequalities are shown to hold for the complex second order elliptic operators \cite{CNW}. We will mainly focus on the inequalities across the interface. Let us first fix some notations.  Assume that $\Omega\subset\R^n$ is an open bounded domain with Lipschitz boundary and $\Sigma\subset\Omega$ is a $C^{1,1}$ hypersurface. Furthermore, assume that $\Omega\setminus\Sigma$ only has two connected components, which we denote $\Omega_\pm$. Let $A_\pm(x)=(a_{\l j}(x)^\pm)_{\l, j=1}^n, W(x), V(x)$ be bounded measurable complex valued coefficients defined in $\Omega$. We say that
\[
\zeta^\pm:=\left(A_\pm, W, V  \right)\in \mathscr{V}(\Omega_\pm,\lambda_0,M_0,K_1, K_2)
\]
if $A_\pm$ satisfy \eqref{symm0}-\eqref{7.4} for $x\in\Omega_\pm$ and $W, V$ satisfy
\[
\|W\|_{L^\infty(\Omega)}\le K_1,\quad \|V\|_{L^\infty(\Omega)}\le K_2.
\]
We will use the notation ${\mathcal L}_\zeta$ to denote
\[
\mathcal{L}_\zeta(x,D)u=\sum_{\pm}H_{\pm}{\rm div}(A_{\pm}(x)\nabla u_{\pm})+\sum_{\pm}H_{\pm}\left(W\cdot\nabla u_{\pm}+Vu_{\pm}\right)\quad\mbox{in}\quad\Omega.
\]
Here, by an abuse of notation, we denote $H_\pm=\chi_{\Omega_\pm}$. 

For an open set $U\subset \R^n$ and a number $s>0$, we define
\begin{equation*}
U^s=\{x\in\R^n:\mbox{dist}(x,U)<s\},
\end{equation*}
\begin{equation*}
U_s=\{x\in U:\mbox{dist}(x,\partial U)>s\},
\end{equation*}
and
\begin{equation*}
sU=\{sx:x\in U\}.
\end{equation*}

\begin{definition}\label{def-lipschitz}
We say that $\dom\in C^{{{k,1}}}$, ${{k\in{\mathbb N}}}$ with constants $\rho_1$, $M_1$ if for any point $P\in\partial\dom$, after a rigid transformation, $P=0$ and 
\begin{equation*}
\dom\cap \Gamma_{\rho_1, M_1}(0)=\{(x',x_n): x'\in \R^{n-1}, |x'|<\rho_1, x_n\in\R, x_n>\Phi(x')\},
\end{equation*}
where $\Phi$ is a $C^{k,1}$ function such that $\Phi(0)=0$,  $\|\Phi\|_{C^{\alpha,1}(B_{\rho_1}(0))}\leq M_1$, and
\begin{equation*}
\Gamma_{\rho_1, M_1}(0)=\{(x',x_n):x'\in \R^{n-1}, |x'|<\rho_1, |x_n|\leq M_1\}.
\end{equation*}
\end{definition}
Throughout this paper, when saying that a domain is $C^{k,1}$, we will mean that it is $C^{k,1}$ with constants $\rho_1$ and $M_1$.

\begin{definition}\label{def-sigma}
We say that $\sig$ is $C^{1,1}$ with constants $\rho_0$, $K_0$ if for any point $P\in\sig$, after a rigid transformation, $P=0$ and 
\begin{equation*}
\dom_\pm\cap C_{\rho_0, K_0}(0)=\{(x',x_n): x'\in \R^{n-1}, |x'|<\rho_0, x_n\in\R, x_n\gtrless\psi(x)\},
\end{equation*}
where $\psi$ is a $C^{1,1}$ function such that $\psi(0)=0$, $\nabla_{x'}\psi(0)=0$,  $\|\psi\|_{C^{1,1}(B_{\rho_0}(0))}\leq K_0$, and
\begin{equation*}
C_{\rho_0,K_0}(0)=\{(x',x_n):x'\in \R^{n-1}, |x'|<\rho_0, |x_n|\leq\frac{1}{2}K_0\rho_0^2\}.
\end{equation*}
\end{definition}

If $\sig$ is as above, then we may "flatten'' the boundary around the point $P$ (without loss of generality $P=0$) via the local $C^{1,1}$-diffeomeorphism
\begin{equation*}
\Psi_P(x,y)= (x, y-\psi(x)).
\end{equation*}
With these assumptions, we will prove a propagation of smallness result as follows.

\begin{theorem}\label{main-thm-1}
Suppose $u\in H^1(\dom)$ solves
\begin{equation*}
{\mathcal L}_\zeta u=0\quad\mbox{in}\quad\Omega.
\end{equation*}
Then there exist $\gamma_0$, depending on $\lambda_0, M_0$, and $h_0$, depending on $\lambda_0$, $M_0$, $K_1$, $K_2$, $\rho_0$, $K_0$, 
such that if $\gamma<\gamma_0$ and $0<h<h_0$,  
 $ r/2>h$, $D\subset\dom$ is connected, open, and {{$D\setminus{\sig}$ has two connected components, denoted by $D_\pm$}}, such that $B_{r}(x_0)\subset D$, $\mbox{\rm dist} (D,\partial\dom)\geq h$, then 
\begin{equation*}
\|u\|_{L^2(D)}\leq C \|u\|_{L^2(B_r(x_0))}^\delta \|u\|_{L^2(\dom)}^{1-\delta},
\end{equation*}
where
\begin{equation*}
C=C_1\left(\frac{|\dom|}{h^{n}}\right)^\frac{1}{2}e^{C_3h^{-s}},\quad \delta\geq\tau^{\frac{C_2|\dom|}{h^{n}}},
\end{equation*}
with $s=s(\lambda_0,K_1, K_2)$, $C_1, C_2>0$, $\tau\in(0,1)$ depending on $\lambda_0$, $M_0$, $K_1$, $K_2$,  $\rho_0$, $K_0$.
\end{theorem}

We would like to remark that the propagation of smallness in Theorem~\ref{main-thm-1} is valid regardless the locations of $D$ and $B_{r_0}(x_0)$, which may intersect the interface $\Sigma$. The strategy of proving Theorem~\ref{main-thm-1} consists two parts.  When we are at one side of the interface, we can use the usual propagation of smallness for equations with Lipschitz complex coefficients based on \cite{CNW}. When near the interface, we then use the three-region inequality derived above to propagate the smallness across the interface. The rest of this section is devoted to the proof of Theorem~\ref{main-thm-1}.

\subsection{Propagation of smallness away from the interface}

In this subsection, we want to derive a propagation of smallness for second order complex elliptic operators with Lipschitz leading coefficients. We consider 
$\zeta:=\left(A, W, V  \right)\in \mathscr{V}(U,\lambda_0,M_0,K_1, K_2)$, where $U$ is an open bounded domain and the value of $\gamma$ in \eqref{complex0} is irrelevant. Note that here $A$ is Lipschitz without jumps in $U$. The following three-ball inequality was proved in \cite[Theorem 3]{CNW}.
\begin{pr}\label{T1}
	Assume that $\zeta:=\left(A, W, V  \right)\in \mathscr{V}(U,\lambda_0,M_0,K_1, K_2)$ and $u\in H_{loc}^1(\Omega)$ solves ${\mathcal L}_\zeta u=0$ in $U$. Then there exist positive constants $R=R(n,\lambda_0,M_0)$  and  $s=s(\lambda_0,K_1, K_2)$ such that if   $0< r_{0} <r_{1}< \lambda_0 r_{2}/2 < \sqrt{\lambda_0}R/2$ with $B_{r_{2}}(x_0)\subset U$, then
	\begin{equation}\label{3ball}
	\|u\|_{L^2(B_{r_1}(x_0))}\leq C\|u\|^\tau_{L^2(B_{r_0}(x_0))}\|u\|^{1-\tau}_{L^2(B_{r_2}(x_0))},
	\end{equation}
	where $C$ is explicitly given by
	\[
	C=e^{C_1(r_0^{-s}-r_2^{-s})}
	\]
	with $C_1>$ depending on $\lambda_0,M_0,K_1, K_2$, and
	\[
	\tau=\frac{(2r_1/\lambda_0)^{-s}-r_2^{-s}}{r_0^{-s}-r_2^{-s}}=\frac{(2r_1/r_2\lambda_0)^{-s}-1}{(r_0/r_2)^{-s}-1}.
	\]
\end{pr}  
The proof of Proposition~\ref{T1} relies on the Carleman estimate derived in \cite{cgt}. 

Having established the three-ball inequality \eqref{3ball}, we can prove the following propagation of smallness based on the chain of balls argument in \cite[Theorem~5.1]{arrv}.  We will not repeat the argument here. 
\begin{pr}\label{T2}
	Assume that $\zeta:=\left(A, W, V  \right)\in \mathscr{V}(U,\lambda_0,M_0,K_1, K_2)$ and $u\in H_{loc}^1(\Omega)$ solves ${\mathcal L}_\zeta u=0$ in $U$. Let $0<h< r/2$ with $r\le\sqrt{\lambda_0}R/2$, $D\subset U$ connected, open, and such that $B_{r}(x_0)\subset D$, $\mbox{\rm dist} (D,\partial U)\geq h$. Then
\begin{equation*}
\|u\|_{L^2(D)}\leq C \|u\|_{L^2(B_r(x_0))}^\delta\|u\|_{L^2(U)}^{1-\delta},
\end{equation*}
where
\begin{equation*}
C=C_2\left(\frac{|U|}{h^n}\right)^{\frac{1}{2}}e^{C_3h^{-s}},\quad \delta\geq\tau^{\frac{C_4|U|}{h^n}},
\end{equation*}
with $C_2, C_3, C_4>0$ depending on $\lambda_0$, $M_0$, $K_1$, $K_2$.
\end{pr}  

\subsection{Propagation of smallness -- an intermediate result}

Here we would like to prove an intermediate propagation of smallness result in which the small ball lies entirely on one side of the interface. Assume that  $D\subset\subset\dom$ is open and connected. {{Recall that we have assumed that $\dom\setminus\sig$ and $D\setminus{\sig}$ both have two connected components, denoted by $\dom_\pm$ and $D_\pm$, respectively.}} Let $\omega_\sig$ be the surface measure induced on $\sig$ by the Lebesgue measure on $\R^n$. We will consider coefficients
\begin{equation*}
\zeta=\left(A_\pm, W,V  \right)\in \mathscr{V}(\dom_\pm,\lambda_0,M_0,K_1, K_2).
\end{equation*}
We can now prove the following propagation of smallness result. 

\begin{theorem}\label{propagation-thm}
Suppose $u\in H^1(\dom)$ solves $\el_{\zeta} u=0$ in $\Omega$. Then there exist $\gamma_0$, depending on $\lambda_0, M_0$, and $h_0$, depending on $\lambda_0$, $M_0$, $K_1$, $K_2$, $\rho_0$, $K_0$, such that if $\gamma<\gamma_0$ and $0<h\le h_0$, $h< r/2$,  $B_{r}(x_0)\subset D_+$, and $\mbox{\rm dist} (D,\partial\dom)\geq h$, then 
\begin{equation*}
\|u\|_{L^2(D)}\leq C \|u\|_{L^2(B_r(x_0))}^\delta\|u\|_{L^2(\dom)}^{1-\delta},
\end{equation*}
where
\begin{equation*}
C=C_1\left(\frac{|\dom|}{h^{n}}\right)\left[1+\left(\frac{\omega_\sig(\sig\cap\dom)}{h^{n-1}} \right)^\frac{1}{2}\right]e^{C_3h^{-s}},\quad \delta\geq\tau^{\frac{C_2|\dom|}{h^{n}}},
\end{equation*}
with $C_1, C_2, C_3>0$, $\tau\in(0,1)$ depending on $\lambda_0$, $M_0$, $K_1$, $K_2$, $\rho_0$, $K_0$.
\end{theorem}

The difficult part of proving Theorem~\ref{propagation-thm} is to obtain $L^2$ estimates of the solution in a neighborhood of $\sig$. We will use Theorem~\ref{thm9.1} to overcome this difficulty. However, we cannot apply Theorem~\ref{thm9.1} directly. The family of regions given in Theorem~\ref{thm9.1} has one serious drawback. If we choose the parameters $R_1=\theta \bar R_1$, $R_2=\theta\bar R_2$, $\theta\in(0,1)$, the vertical sizes of the regions would scale like $\theta$, while their horizontal sizes would scale like $\theta^{\frac{1}{2}}$. Using just these two parameters in the proof would then lead to constants in the propagation of smallness inequality (i.e. the constants $C$ and $\delta$ in Theorem \ref{main-thm-1}) that depend on the geometry of $\dom$, $D$, and $B_{r}(x_0)$ in a way that is not invariant under a rescaling of these sets. Therefore, we will study how the three-region inequality \eqref{9.1} behaves under scaling. 

Let us first introduce the scaled coefficients
\begin{equation*}
\tilde\zeta=\left(\tilde A_\pm, \tilde W, \tilde V  \right)\in \mathscr{V}(\R^n_\pm,\lambda_0,M_0,K_1, K_2).
\end{equation*}
For $0<\theta\leq 1$, let
\[
\mathcal{L}^\theta_{\tilde\zeta}(\cdot,D)v=\sum_{\pm}H_{\pm}{\rm div}(\tilde A_{\pm}(\theta \cdot)\nabla v_{\pm})+\sum_{\pm}H_{\pm}\left(\theta\tilde W(\theta\cdot)\cdot\nabla v_{\pm}+\theta^2\tilde V(\theta\cdot)v_{\pm}\right).
\]
Note that if $\tilde\zeta=\left(\tilde A_\pm, \tilde W, \tilde V  \right)\in \mathscr{V}(\R^n_\pm,\lambda_0,M_0,K_1, K_2)$, then
\begin{equation*}
\left(\tilde A_\pm(\theta\cdot), \theta\tilde W(\theta\cdot), \theta^2\tilde V(\theta\cdot)  \right)\in\mathscr{V}(\R^{n}_\pm, \lambda_0,\theta M_0,\theta K_1, \theta^2 K_2).
\end{equation*}
It is also clear that if $\el_{\tilde\zeta}u=0$ in $\dom$, then $u^\theta(x)=\theta^{-2}u(\theta x)$ solves
\begin{equation*}
\mathcal{L}^\theta_{\tilde\zeta} u^\theta=0\quad\mbox{in}\quad\theta^{-1}\dom. 
\end{equation*}
Moreover, if $U\subset \theta^{-1}\dom$, we have
\begin{equation*}
\int_{\theta U}|u(x)|^2dx=\theta^{n+4}\int_U |u^\theta(y)|^2dy.
\end{equation*}
We therefore obtain, by scaling, the following corollary to Theorem \ref{thm9.1}.
\begin{pr}\label{thm-FLVW-2}
Assume that the assumptions in Theorem~\ref{thm9.1} hold. Let $0<R_1, R_2\leq R$, $\theta\in (0,1]$, and 
\begin{equation*}
\el_{\tilde\gamma} u=0\quad\mbox{in}\quad \theta U_3,
\end{equation*}
then 
\begin{equation*}
\int_{\theta U_2}|u|^2\leq (e^{\tau_0 R_2}+CR_1^{-4})
\left(\int_{\theta U_1} |u|^2\right)^{\frac{R_2}{2R_1+3R_2}}
\left(\int_{\theta U_3} |u|^2\right)^{\frac{2R_1+2R_2}{2R_1+3R_2}}.
\end{equation*}
\end{pr}

In order to adapt that result to the possibly curved surface $\sig$, we need to first consider how the three regions transform under a local boundary flattening diffeomorphism $\Psi_P$. Pick a point $P\in\sig$ and set $P=0$ without loss of generality. Let  $(x',x_n)\in C_{\rho_0,K_0}(0)$. We will try to determine when $(x',x_n)\in\Psi^{-1}_P(\theta U_2)$. To this end, we introduce the notation
\begin{equation*}
y'=x',\quad y_n=x_n-\psi(x').
\end{equation*}
It is clear that $(x',x_n)\in\Psi^{-1}_P(\theta U_2)$ if and only if $\theta^{-1}(y',y_n)\in U_2$. We denote
\[
\eta(x')=\frac{\psi(x')}{|x'|^2},
\]
which is a bounded function due to the regularity assumption of $\Sigma$. It was proved in \cite[Lemma~3.1]{CW} that if
\begin{equation}\label{u2}
r<\theta\min\left\{\frac{\delta R_1}{6a\alpha_-},\frac{2\delta R_2}{3\alpha_-}, \frac{R_1}{12a},\theta^{-1}\rho_0,\rho_1,\rho_2,\rho_3\right\},
\end{equation}
then $\Psi_P(B_r(P))\subset \theta U_2$, where
\[
\left\{
\begin{aligned}
&\rho_1=\frac{\alpha_-\delta}{\delta+\beta},\\
&\rho_2\;\mbox{is chosen such that}\; 2\|\eta\|\rho_2+\|\eta\|^2\rho_2^2<\frac 12,\; \|\eta\|=\|\eta\|_{L^\infty(B'_{\rho_0}(0))},\\
&\rho_3=\frac{2\alpha_-\delta}{\beta}.
\end{aligned}\right.
\]
In \cite[Lemma~3.2]{CW}, the following relation was established. 
\begin{lemma}\label{lem-U3}
$\Psi_P^{-1}(\theta U_3)$ is contained in a ball of radius
{\scriptsize
\begin{equation*}
\theta\left[(1+2\|\eta\|^2)\left(\frac{2\alpha_-}{a}R_1+8\delta R_2\right)+
\frac{1}{a^2}\left[2+(1+2\|\eta\|^2)\frac{\beta}{\delta}   \right]R_1^2 +
\frac{128\delta^2R_2^2}{\left[\alpha_-+\sqrt{\alpha_-^2-8\beta  R_2}  \right]^2}      \right]^{1/2}
\end{equation*}
}
centered at $P$.
\end{lemma}

Finally, we need to estimate the distance from $\Psi_P^{-1}(\theta U_1)$ to $\sig\cap C_{\rho_0,K_0}$. It was proved in \cite[Lemma~3.3]{CW} that
\begin{equation}\label{u1}
\mbox{\rm dist}(\Psi_P^{-1}(\theta U_1),\Sigma)>\theta\frac{R_1}{16a}.
\end{equation}

We are now ready to prove Theorem \ref{propagation-thm}. We will follow the arguments used in the proof of Theorem~3.1 in \cite{CW}. Since we have slightly different constants here, we provide the proof for the sake of completeness.

\medskip\noindent
{\bf Proof of Theorem}~\ref{propagation-thm}. By the assumption, we may take $D$ to be the set 
\begin{equation*}
D=\{x\in\dom:\mbox{\rm dist}(x,\partial\dom)>h\}.
\end{equation*}
We want to point out that even though the choice of $\alpha_\pm$ in Theorem~\ref{thm8.2} depends on $A_\pm(P)$ for $P\in\Sigma$, we can choose a pair of $\alpha_\pm$ such that Carleman estimate \eqref{8.24} holds near all $P\in\Sigma$ in view of the regularity assumptions of $A_\pm$ and $\Sigma$. Consequently, we can pick $R_1$, $R_2$ so that we can apply Proposition~\ref{thm-FLVW-2} at any point $P\in \sig\cap D$. By Lemma \ref{lem-U3}, there is a constant $d>0$, independent of $P$, such that $\Psi_P^{-1}(\theta U_3)\subset B_{\theta d}(P)$. We then choose $\theta$ such that $\theta d=\frac{h}{2}$, which implies $\Psi_P^{-1}(\theta U_3)\subset \dom$ for any $P\in \sig\cap D$. Of course, this choice is not possible if $h$ is too large. Therefore, we need to set $h_0$ small enough, depending on $\rho_0, K_0, \lambda_0, M_0, K_1, K_2$.

With this choice of parameters, by \eqref{u1}, there is a constant $0<\mu<1$, also independent on $P$, so that 
\begin{equation*}
\mbox{\rm dist}(\Psi_P^{-1}(\theta U_1), \sig)>\mu h.
\end{equation*}
 Note that, depending on the geometry of $\sig$, we again need to set $h_0$ and $R$ small enough so that $\Psi_P^{-1}(\theta U_1)\cap \sig^{\mu h}=\emptyset$, for any  $P\in \sig\cap D$.
 
It follows from \eqref{u2} that there exists a constant $\nu>0$, and without loss of generality $\nu<\mu<1$, such that $B_{5\nu h}(P)\subset \Psi_P^{-1}(\theta U_2)$. By Vitali's covering lemma, there exist finitely many $P_1,\ldots,P_N\in \sig\cap D$ so that
\begin{equation}\label{qq}
\sig^{\nu h}\cap D\subset\bigcup_{j=1}^N\Psi_{P_j}^{-1}(\theta U_2),
\end{equation}
and the balls $B_{\nu h}(P_j)$ are pairwise disjoint. By this last property, since for small $h$ we have $\omega_\sig(\sig^{\nu h}\cap D)\sim \nu h \omega_\sig(\sig\cap D)$, it follows that there is a constant $C$ such that
\begin{equation}\label{qqq}
N\leq C\frac{\omega_\sig(\sig\cap D)}{h^{n-1}}\leq C\frac{\omega_\sig(\sig\cap\dom)}{h^{n-1}}.
\end{equation}

Let us denote $\tilde D= (D_+)^{h/2}\setminus(\sig^{\nu h}\cup \dom_-)$, then by Proposition \ref{T2}, we have that
\begin{equation}\label{D-plus-estimate}
\|u\|_{L^2(\tilde D)}\leq C_+ \|u\|_{L^2(B_r(x_0))}^{\delta_+}\|u\|_{L^2(\dom)}^{1-\delta_+},
\end{equation}
where
\begin{equation*}
B_r(x_0)\subset\tilde D,\quad C_+=C_2\left(\frac{|\dom|}{h^n}\right)^{\frac{1}{2}}e^{C_3h^{-s}},\quad \delta_+\geq\tau^{\frac{C_4|\dom|}{h^n}}.
\end{equation*}
The function $v=u\circ \Psi_{P_j}^{-1}\in H^1(\theta U_3)$ satisfies $\el_{\tilde \zeta} v=0$ in $\theta U_3$ with $\tilde \zeta$ satisfying
\begin{equation*}
\tilde \zeta^\pm=\left(\tilde A_\pm,  \tilde W, \tilde V  \right)\in \mathscr{V}(\R^n_\pm,\tilde \lambda_0,\tilde M_0,\tilde K_1, \tilde K_2),
\end{equation*}
with $C>0$ and the parameters $\tilde \lambda_0,\tilde M_0,\tilde K_1,\tilde K_2$ depending on $\lambda_0, M_0, K_1, K_2, \rho_0, K_0$. We can then pull back the three regions inequality in Proposition \ref{thm-FLVW-2} and apply it to $u$ and the regions $\Psi_{P_j}^{-1}(\theta U_1)$, $\Psi_{P_j}^{-1}(\theta U_2)$, $\Psi_{P_j}^{-1}(\theta U_3)$.

Since 
\begin{equation*}
\Psi_{P_j}^{-1}(\theta U_1)\subset (D_+)^{h/2}\setminus(\sig^{\nu h}\cup \dom_-),
\end{equation*}
we have that
\begin{equation*}
\|u\|_{L^2(\Psi_{P_j}^{-1}(\theta U_2))}\leq C \|u\|_{L^2((D_+)^{h/2}\setminus(\sig^{\nu h}\cup \dom_-))}^\xi 
\|u\|_{L^2(\dom)}^{1-\xi},
\end{equation*}
where $\xi=\frac{R_2}{2 R_1+3 R_2}$. Combining this and \eqref{D-plus-estimate}, we obtain
\begin{equation*}
\|u\|_{L^2(\Psi_{P_j}^{-1}(\theta U_2))}\leq C_1'\left(\frac{|\dom|}{h^n}\right)^{\frac{\xi}{2}}e^{C_3\xi h^{-s}}
 \|u\|_{L^2(B_r(x_0))}^{\xi\delta_+}\|u\|_{L^2(\dom)}^{1-\xi\delta_+}.
\end{equation*}
Then it follows from \eqref{qq} and \eqref{qqq} that 
\begin{equation}\label{Sigma-estimate}
\|u\|_{L^2(\sig^{\nu h}\cap D)}\leq C_1''\left(\frac{\omega_\sig(\sig\cap\dom)}{h^{n-1}}\right)^{\frac{1}{2}}\left(\frac{|\dom|}{h^n}\right)^{\frac{\xi}{2}}e^{C_3\xi h^{-s}}\|u\|_{L^2(B_r(x_0))}^{\xi\delta_+}\|u\|_{L^2(\dom)}^{1-\xi\delta_+}.
\end{equation}

Applying Proposition \ref{T2} again (now with an appropriate small ball $\tilde B_{\tilde r}\subset\Sigma^{\nu h}\cap D_-\subset\Sigma^{\nu h}\cap D$), we have
\begin{equation}\label{D-minus-estimate}
\begin{aligned}
\|u\|_{L^2(D_-\setminus\sig^{\nu h})}\leq &C_1'''\left(\frac{|\dom|}{h^n}  \right)^{\frac{1}{2}}\left(\frac{\omega_\sig(\sig\cap\dom)}{h^{n-1}}\right)^{\frac{\delta_-}{2}}\left(\frac{|\dom|}{h^n}\right)^{\frac{\delta_-\xi}{2}}e^{C_3(1+\delta_-\xi) h^{-s}}\\
&\times \|u\|_{L^2(B_r(x_0))}^{\delta_-\xi\delta_+}\|u\|_{L^2(\dom)}^{1-\delta_-\xi\delta_+},
\end{aligned}
\end{equation}
where
\begin{equation*}
\delta_-\geq\tau^{\frac{C_2'|\dom|}{h^n}}.
\end{equation*}
Combining estimates \eqref{D-plus-estimate}, \eqref{Sigma-estimate}, and \eqref{D-minus-estimate}, we obtain the conclusion of Theorem \ref{propagation-thm}.
\eproof

\subsection{Propagation of smallness -- a general result}

This section is devoted to the proof of Theorem~\ref{main-thm-1}. We will first prove a three-ball inequality, which is a direct consequence of Theorem \ref{propagation-thm}. Such three balls inequality is a building block in the proof of Theorem~\ref{main-thm-1}. 

\subsubsection{Three-ball inequality}

Recall that we  assume $\dom\subset\R^n$ is an open Lipschtiz domain,  $\sig$ is a $C^{1,1}$ hypersurface with constants $\rho_0$, $K_0$,  and $\dom\setminus \sig$ has two connected components, $\dom_\pm$. We consider coefficients
\begin{equation*}
\zeta=\left(A_\pm,  W, V  \right)\in \mathscr{V}(\dom_\pm,\lambda_0,M_0,K_1, K_2)
\end{equation*}
with $\gamma<\gamma_0$. Also, let $u\in H^1(\dom)$ be a solution of
\begin{equation*}
\el_{\zeta} u=0\quad\mbox{in}\quad\Omega.
\end{equation*}

\begin{pr}\label{three-balls}
There exists $\bar r>0$, depending on $\rho_0$, $K_0$, such that if  $0<r_1<r_2<r_3<\bar r$, $Q\in\dom$, $\text{\rm dist}(Q,\partial\dom)>r_3$, then there exist $C>0$, $0<\delta<1$ such that
\begin{equation}\label{3-ball}
\|u\|_{L^2(B_{r_2}(Q))}\leq C\|u\|_{L^2(B_{r_1}(Q))}^\delta \|u\|_{L^2(B_{r_3}(Q))}^{1-\delta}
\end{equation}
with 
\[
C=C_1e^{C_2r_3^{-s}},
\]
where $C_1, C_2$, and $\delta$ depend on  $\lambda_0$, $M_0$, $\rho_0$, $K_0$, $K_1$, $K_2$, $\frac{r_1}{r_3}$, $\frac{r_2}{r_3}$, $\text{\rm diam}(\dom)$. 
\end{pr}

\pf.
We would like to use the propagation of smallness result in Theorem \ref{propagation-thm} with $r=\frac{r_1}{10}$, $D=B_{r_2}(Q)$, and replacing $\dom$ there by $B_{r_3}(Q)$. We can choose the constant $\bar r$ so  that $B_{r_j}(Q)\setminus\sig$ can all only have at most two connected components. This would be the case for example if $\bar r\leq \min\{\rho_0,\frac{1}{2}K_0\rho_0^2\}$. Fix an $\bar r$ as described. Then we can always find $Q'\in B_{r_1}(Q)$ so that $B_{r_1/10}(Q')\subset B_{r_1}(Q)\cap\dom_+$ or $B_{r_1/10}(Q')\subset B_{r_1}(Q)\cap\dom_-$. Without loss of generality we may assume that $B_{r_1/10}(Q')\subset B_{r_1}(Q)\cap\dom_+$.

Let $g^\sig$ be the metric induced on $\sig$ by the Euclidean metric of $\R^n$. Around a point $P\in\sig$ at which we have chosen coordinates as in Definition \ref{def-sigma}, we can use the coordinates $(x_1, \ldots, x_{n-1})$ as a local map for $\sig$. In these coordinates, we have
\begin{equation*}
g^\sig_{jk}=\delta_{jk}+\partial_j\psi\partial_k\psi.
\end{equation*}
This observation implies that
there exists a constant $\kappa$ so that
\begin{equation*}
\omega_\sig(\sig\cap B_{r_3}(Q))<\kappa r_3^{n-1}.
\end{equation*}

Here we can also observe that we can cover $\sig\cap\dom$ with balls $B_{\rho_0}(P)$, $P\in \dom$, and on each of those the same estimate for $\omega_\sig(\sig\cap B_{\rho_0}(P))$ holds. Using again Vitali's covering lemma, it is not hard to see that we can make sure that the number $N$ of these balls satisfies 
\begin{equation*}
N<C\,\text{\rm diam}(\dom)^n,
\end{equation*} 
with $C$ depending only on the dimension $n$. It follows that $\omega_\sig(\sig\cap\dom)$ can be bound from above by a quantity depending on $\rho_0$, $K_0$, and $\text{\rm diam}(\dom)$.

We will treat several cases separately. The first case is when $r_3-r_2<\min\{\frac{r_1}{20},h_0\}$. Remind that the constant $h_0$ is given in Theorem \ref{propagation-thm}. We now apply Theorem \ref{propagation-thm}, with $h=r_3-r_2$,  to obtain
\begin{equation*}
\|u\|_{L^2(B_{r_2}(Q))}\leq C\|u\|_{L^2(B_{\frac{r_1}{10}}(Q'))}^\delta \|u\|_{L^2(B_{r_3}(Q))}^{1-\delta},
\end{equation*}
where
\begin{equation}\label{C-i}
C=C_1\left(\frac{r_3^n}{(r_3-r_2)^{n}}\right)\left[1+\left(\frac{\kappa r_3^{n-1}}{(r_3-r_2)^{n-1}} \right)^\frac{1}{2}\right]e^{C_3(r_3-r_2)^{-s}},
\end{equation}
\begin{equation}\label{delta-i}
 \delta\geq\tau^{C_2\frac{ r_3^n}{(r_3-r_2)^{n}}}.
\end{equation}

The second case is when $\frac{r_1}{10}<2h_0$, $r_3-r_2\geq\frac{r_1}{20}$. Let $r_3'=r_2+\frac{r_1}{21}$ (note $r_3'<r_3$), $h=\frac{r_1}{21}$ and again apply Theorem \ref{propagation-thm} to obtain
\begin{equation*}
\|u\|_{L^2(B_{r_2}(Q))}\leq C\|u\|_{L^2(B_{\frac{r_1}{10}}(Q'))}^\delta \|u\|_{L^2(B_{r_3}(Q))}^{1-\delta},
\end{equation*}
where
\begin{equation}\label{C-ii}
C=C_1\left(\frac{(r_2+r_1/21)^n}{(r_1/21)^{n}}\right)\left[1+\left(\frac{\kappa (r_2+r_1/21)^{n-1}}{(r_1/21)^{n-1}} \right)^\frac{1}{2}\right]e^{C_3(r_1/21)^{-s}},
\end{equation}
\begin{equation}\label{delta-ii}
 \delta\geq\tau^{C_2\frac{ (r_2+r_1/21)^n}{(r_1/21)^{n}}}.
\end{equation}

The third and final case is when $\frac{r_1}{10}\geq 2h_0$, $r_3-r_2\geq h_0$. In this case we take $h=h_0$, and use the estimates
\begin{equation*}
|B_{r_3}(Q)|\leq (\mbox{\rm diam}(\dom))^n,\quad \omega_\sig(B_{r_3}(Q)\cap\sig)\leq \omega_\sig(\sig\cap\dom).
\end{equation*}
We then have
\begin{equation*}
||u||_{L^2(B_{r_2}(Q))}\leq C\|u\|_{L^2(B_{\frac{r_1}{10}}(Q'))}^\delta \|u\|_{L^2(B_{r_3}(Q))}^{1-\delta},
\end{equation*}
where
\begin{equation}\label{C-iii}
C=C_1\frac{(\mbox{\rm diam}(\dom))^n}{h_0^{n}}\left[1+\left(\frac{ \omega_\sig(\sig\cap\dom)}{h_0^{n-1}} \right)^\frac{1}{2}\right]e^{C_3h_0^{-s}},
\end{equation}
\begin{equation}\label{delta-iii}
 \delta\geq\tau^{C_2\frac{ (\text{diam}(\dom))^n}{h_0^{n}}}.
\end{equation}
It follows that, in all cases, we have our three ball inequality with the constant $C$ being the maximum of the ones in \eqref{C-i}, \eqref{C-ii}, and \eqref{C-iii}, and the exponent $\delta$ being the minimum of the ones in \eqref{delta-i}, \eqref{delta-ii}, and \eqref{delta-iii}.

\subsubsection{Proof of Theorem \ref{main-thm-1}}

Once we have established the three balls inequality in Theorem \ref{three-balls}, the proof of Theorem \ref{main-thm-1} is standard. We include it here for the benefit of the reader. Let
\begin{equation*}
r_3=\frac{h}{2},\quad r_2=\frac{1}{5}r_3=\frac{1}{10}h,\quad r_1=\frac{1}{3}r_3=\frac{1}{30}h,
\end{equation*}
and
\begin{equation*}
\tilde D=\left\{x\in\dom:  \mbox{\rm dist}(x,D)<r_1\right\},
\end{equation*}
which is an open connected subset of $\dom$, such that $D\subset \tilde D$, $\text{\rm dist}(\tilde D,\partial\dom)>h/2$. Let $y\in\tilde D$ and $c(t)\in C([0,1];\tilde D)$ be a continuous curve such that $c(0)=x_0$, and $c(1)=y$. Define 
\begin{equation*}
0=t_0<t_1<\cdots<t_N=1
\end{equation*}
so that
\begin{equation*}
\begin{array}{c}
t_{k+1}=\max\{t:|c(t)-c(t_k)|=2r_1\}, \text{ as long as } |y-c(t_k)|>2r_1,\\[5pt] \text{ otherwise }N=k+1, t_N=1.
\end{array}
\end{equation*}
Then $B_{r_1}(c(t_k))\cap B_{r_1}(c(t_{k-1}))=\emptyset$, and $B_{r_1}(c(t_{k+1}))\subset B_{r_2}(c(t_{k}))$, $k=1,\ldots,N-1$. By Theorem \ref{three-balls} we have
\begin{equation*}
\|u\|_{L^2(B_{r_1}(c(t_{k+1})))}\leq C\|u\|_{L^2(B_{r_1}(c(t_{k})))}^\tau\|u\|_{L^2(\dom)}^{1-\tau},
\end{equation*}
where $k=0,\ldots,N-1$.  Note that $C=C_1e^{C_2h^{-s}}$ with $C_1,C_2$ depending on $\lambda_0$, $M_0$, $\rho_0$, $K_0$, $K_1$, $K_2$, $\text{\rm diam}(\dom)$.

Let
\begin{equation*}
m_k=\frac{\|u\|_{L^2(B_{r_1}(c(t_k))}}{\|u\|_{L^2(\dom)}},
\end{equation*}
then $m_{k+1}\leq C m_k^\tau$, $k=0,\ldots,N-1$, and so
\begin{equation*}
m_N\leq C^{1+\tau+\cdots+\tau^{N-1}}m_0^{\tau^N}.
\end{equation*}
Since the balls $B_{r_1}(c(t_k))$ are pairwise disjoint,
\begin{equation*}
N\leq \frac{|\dom|}{\omega_n r_1^n}\leq \frac{C_2|\dom|}{h^n}.
\end{equation*}
Then it is easy to see that 
\begin{equation*}
\tau^N\geq\tau^{\frac{C_2|\dom|}{h^n}},\quad C^{1+\tau+\cdots+\tau^{N-1}}\leq C^{\frac{1}{1-\tau}}.
\end{equation*}
From a family of disjoint open cubes of side $2r_1/\sqrt{n}$ whose closures cover $\R^n$, extract the finite number of cubes which intersect $D$ non-trivially: $Q_j$, $j=1,\ldots,J$. The number of these cubes satisfies $J\leq\frac{n^{n/2}|\dom|}{2^nr_1^n}$. For each $j$ there exists $w_j\in \tilde D$ such that $Q_j\subset B_{r_1}(w_j)$. Then
\begin{equation*}
\int_D |u|^2\leq \sum_{j=1}^J\int_{Q_j}|u|^2\leq
\sum_{j=1}^J\int_{B_{r_1}(w_j)}|u|^2\leq
JC^{2/(1-\tau)}\|u\|_{L^2(B_{r_1}(x_0))}^{2\delta}\|u\|_{L^2(\dom)}^{2(1-\delta)}.
\end{equation*}
\eproof

\section{Size estimate}\label{size-est}

In this section, we will study the size estimate problem using the uniqueness estimate derived above. Precisely speaking, we will extend the result obtained for the Lipschitz complex coefficients in \cite{CNW} to the piecewise Lipschitz complex coefficients. So we follow closely the arguments used in \cite{CNW}. 

Let $\Omega$ be a conducting body with an anisotropic background with current-voltage relation (or Ohm's law)  given by
\begin{equation}\label{bg}
I_0(\nabla u)=(\sigma_0(x)+i\eps_0(x))\nabla u(x),
\end{equation}
where $\sigma_0$, $\eps_0$ are symmetric matrix valued functions. As above, we assume that $\Omega\subset\R^n$ is an open bounded domain with Lipschitz boundary with constants $\rho_1, M_1$ and $\Sigma\subset\Omega$ is a $C^{1,1}$ hypersurface with constants $\rho_0, K_0$ such that $\Omega\setminus\Sigma$ only has two connected components, which are denoted by $\Omega_\pm$. Here $\Omega_-$ represents the inner component, i.e., $\partial\Omega_-=\Sigma$. Assume that $\mbox{\rm dist}(\Sigma,\partial\Omega)>0$.

Let $\sigma_0(x)+i\eps_0(x)=A^0_\pm(x)$ be bounded measurable complex valued and 
\[
A^0_\pm\in \mathscr{V}(\Omega_\pm,\lambda_0,M_0),
\]
i.e., $A^0_\pm$ satisfy \eqref{7.2}-\eqref{7.4} for $x\in\Omega_\pm$. Here $\gamma$ in \eqref{complex0} has been chosen so that Theorem~\ref{main-thm-1} holds. 

Let $D\subset \Omega$ represent the region of anomaly hidden in $\Omega$. In \cite{flvw}, $D$ is assumed to be contained entirely in $\Omega_-$, i.e., $\bar D\subset\Omega_-$. It should be emphasized that here the location of $D$ is arbitrary as long as it stays away from $\partial\Omega$ at a fixed distance. It can intersect the interface $\Sigma$. The underlying reason is that here we have more powerful propagation of smallness result, Theorem~\ref{main-thm-1}. Assume that the current-voltage relation in $D$ is given by 
\begin{equation}\label{exotic}
I_1(\nabla u)=(\sigma_1(x)+i\eps_1(x))\nabla u(x)+\zeta_1(x) \overline{\nabla u}(x),
\end{equation}
where $\sigma_1,\eps_1$, $\zeta_1$ are real symmetric matrix-valued functions such that the supports of $\sigma_1-\sigma_0$, $\eps_1-\eps_0$, $\zeta_1$ are contained in $D$. We want to remark that \eqref{exotic} can be derived from the Maxwell equations describing the propagation of electromagnetic waves traveling in a chiral medium. We refer the reader to \cite[Section 1]{CNW} for detailed explanation. Here $\zeta_1$ denotes the chirality of the medium.

Given a Neumann boundary data $g\in L^2(\partial\Omega)$ satisfying $\int_{\partial\Omega}g=0$,  let $u_0$ be the unique solution of the unperturbed equation 
\begin{equation}\label{unpert}
\left\{
\begin{aligned}
&\nabla\cdot I_0(\nabla u_0)=0\quad\text{in}\quad\Omega,\\
&I_0(\nabla u_0)\cdot\nu=g\quad\text{on}\quad\partial\Omega,\\
&\int_{\Omega}u_0=0,
\end{aligned}
\right.
\end{equation}
and $u_1$ be the solution of the perturbed equation with the same boundary data, i.e.
\begin{equation}\label{pert}
\left\{
\begin{aligned}
&\nabla\cdot I_1(\nabla u_1)=0\quad\text{in}\quad\Omega,\\
&I_1(\nabla u_1)\cdot\nu=g\quad\text{on}\quad\partial\Omega,\\
&\int_{\Omega}u_1=0.
\end{aligned}
\right.
\end{equation}
The inverse problem studied here is to estimate the size of the inclusion $D$ by one pair of boundary measurement $\{u|_{\partial\Omega},g\}$. More precisely, one tries to estimate  $|D|$ using the power
gap $\delta W=W_0-W_1$, where
\[
W_1=\int_{\partial\Omega}u_1 g,\quad W_0=\int_{\partial\Omega}u_0 g.
\]

Before stating the main result of the size estimate, we list several conditions imposed on $\sigma_1,\eps_1,\zeta_1$. We denote $Id$ the $n$-dimensional identity matrix.
\begin{itemize}
	\item (Boundedness and Ellipticity) for some $\lambda_1\le 1$,
	\begin{equation}\label{se0}
	\begin{aligned}
	 &\lambda_1 Id\le \sigma_ 1 + \zeta_1\le\lambda_1^{-1} Id,\quad\lambda_1 Id\le \sigma_ 1 - \zeta_1\le\lambda_1^{-1} Id,\\
	 &\|\eps_0\|_{L^\infty(\Omega)}\le\lambda_1^{-1}, \quad \|\eps_1\|_{L^\infty(\Omega)}\le\lambda_1^{-1}.
	\end{aligned}  
	\end{equation}
	\item (Jump condition) There exists $\varrho>0$ such that either
	\begin{equation}\label{a0}
	\begin{aligned}
	 &{\rm (i)} \qquad\zeta_1 \le (\sigma_1-\sigma_0)-\varrho Id,\;\;\zeta_1 \le (\sigma_0-\sigma_1)-\varrho Id \quad\mbox{a.e. } \, x\in D, \\
	 & \mkern-36mu \mbox{or} \\
	 &{\rm (ii)} \qquad \zeta_1 \ge (\sigma_1-\sigma_0)+\varrho Id,\;\;\zeta_1 \ge (\sigma_0-\sigma_1)+\varrho Id \quad\mbox{a.e. } \, x\in D.
	\end{aligned}
	\end{equation}
	Morever, for $\delta=\delta(\alpha,\varrho)>0$ sufficiently small, we have 
	\begin{equation}\label{delta}
	\|\eps_1-\eps_0\|_{L^\infty(D)}\le\delta.
	\end{equation}
\end{itemize}
Note that $\|\eps_0\|_{L^\infty(\Omega)}\le\lambda_1^{-1}$ provided $\gamma\lambda_0^{-1}\le\lambda_1^{-1}$. Let us remark that these assumptions are rather standard in this area, see e.g. \cite{ars} and \cite{bfv}. The main difference between ours and those in  previous work is the jump condition which is adapted to our physical model \eqref{exotic}. It plays an important role in the derivation of the power gap (see Proposition \ref{energy}).

We are ready to state the main theorem about the size estimates. Here we define the free energy
\[
W'_0=\int_{\partial\Omega}\bar u_0g=\int_\Omega\sigma_0\nabla u_0\cdot\nabla\bar u_0+i\int_\Omega\eps_0\nabla u_0\cdot\nabla\bar u_0.
\]
Since $u_0$ is completely determined, $W'_0$ is known.
\begin{theorem}\label{size}
Let $\Omega$ be a bounded domain in $\mathbb{R}^n$ with $\partial\Omega\in C^{3,1}$. Assume that for some $d_0,d_1>0$, the inclusion  $D\subset\Omega$  satisfies $\text{\rm dist}(D,\partial\Omega)\ge d_0$ and 
\begin{equation}\label{i1}
\quad|D_{d_1}|\ge\frac 12|D|\quad\mbox{\rm (fatness condition)}.
\end{equation}
Suppose that \eqref{se0}-\eqref{delta} are satisfied. Then there exist two positive constants $C_1$, and $C_2$ such that
\[
C_1\frac{|\Re\delta W|}{\Re W'_0}\le|D|\le C_2\frac{|\Re\delta W|}{\Re W'_0},
\]
where $C_1$ depends on $d_0,\lambda_0, M_0, \lambda_1, \rho_0, K_0 $ and $C_2$ depends on $\lambda_0, M_0$, $\lambda_1$, $\rho_0, K_0$, $\rho_1$, $K_1$, $\varrho$, $|\Omega|$, $d_0$, $d_1$, and $\|g\|_{L^2(\partial\Omega)}/\|g\|_{H^{-1/2}(\partial\Omega)}$.
\end{theorem}

We remark that  the fatness condition \eqref{i1} is only used to obtain an upper bound of $|D|$. In \cite{bfv}, an estimate was obtained for general inclusions without the fatness condition \eqref{i1}. This was possible because under their conditions, $u_0$ satisfies a doubling inequality, which is not known in the case considered here. The proof of Theorem~\ref{size} relies on the following energy inequalities proved in \cite{CNW}. We include the proof here for readers' convenience.  
\begin{pr}{\rm\cite[Proposition 5]{CNW}}\label{energy}
Assume that \eqref{se0} and \eqref{delta} hold. Then there exist positive constants $\kappa_1=\kappa_1(\lambda_0,\lambda_1,\varrho)$ and $\kappa_2=\kappa_2(\lambda_0,\lambda_1) $ such that
\begin{itemize}
	\item If \eqref{a0}{\rm (i)} holds then
	\begin{equation}\label{eineq1}
	\kappa_1\int_D|\nabla u_0|^2\le \Re\delta W \le \kappa_2\int_D|\nabla u_0|^2.
	\end{equation}
	\item If \eqref{a0}{\rm (ii)} holds then
	\begin{equation}\label{eineq2}
	\kappa_1\int_D|\nabla u_0|^2\le -\Re\delta W \le \kappa_2\int_D|\nabla u_0|^2.
	\end{equation}
\end{itemize}
\end{pr}
\pf.
	In this proof, $C$ denotes a positive constant whose value may change from line to line. We have, for $j=0,1$,
	\begin{equation}\label{coneq1}
	\begin{pmatrix} \Re I_j(\nabla u_j) \\ \Im I_j(\nabla u_j) \end{pmatrix} = \begin{pmatrix}\sigma_j +\zeta _j &-\eps_j\\ \eps_j&\sigma_j-\zeta _j\end{pmatrix} \begin{pmatrix}\Re \nabla u_j\\ \Im \nabla u_j\end{pmatrix},
	\end{equation}
	where we set $\zeta_0=0$. Following Cherkaev and Gibiansky \cite{cg94},  \eqref{coneq1} can be expressed as
	\begin{equation}\label{coneq2}
	\begin{pmatrix} \Re \nabla u_j\\ \Im I_j(\nabla u_j) \end{pmatrix}=\begin{pmatrix}(\sigma_j+\zeta_j)^{-1}&(\sigma_j+\zeta_j)^{-1}\eps_j\\ \eps_j(\sigma_j+\zeta_j)^{-1}&\sigma_j-\zeta_j+\eps_j(\sigma_j+\zeta_j)^{-1}\eps_j\end{pmatrix}\begin{pmatrix}\Re I_j(\nabla u_j) \\ \Im \nabla u_j\end{pmatrix}.
	\end{equation}
	We denote the square matrix on the right hand side $B_j$ and $v_j=(\Re{I}_j(\nabla u_j),\Im \nabla u_j)^t$.  From our assumptions, it can be easily checked that $B_j$ is positive-definite. Using \eqref{coneq2}, we have, for $j,k\in\{0,1\}$,
	\begin{equation}\label{basic}
	\begin{aligned}
	\int_\Omega B_j v_j\cdot v_k &=  \int_\Omega \Re\nabla u_j \cdot \Re{I}_k(\nabla u_k) + \Im{I}_j(\nabla u_j) \cdot \Im \nabla u_k\\
	&=  \Re\int_\Omega \nabla\cdot({I}_k(\nabla u_k) \Re u_j) + \Im\int_\Omega \nabla\cdot({I}_j(\nabla u_j) \Im u_k)\\
	&=\int _{\partial\Omega} ( \Re u_j \Re h + \Im u_k \Im h).
	\end{aligned}
	\end{equation}
	From this and $B_1=B_1^T$, we can deduce that
	\begin{align}\label{basic*}
	\Re \delta W= \Re\int _{\partial\Omega} (\bar u_0-\bar u_1)h =\int_\Omega (B_0-B_1) v_0\cdot v_1= \int_D (B_0-B_1) v_0\cdot v_1.  
	\end{align}
	Using \eqref{basic} and $B_0=B_0^T$, we obtain 
	\begin{equation}
	\begin{aligned}
	\int_\Omega B_0 (v_0-v_1)\cdot(v_0-v_1)  =& \int_\Omega ( B_0 v_0\cdot v_0 -2 B_0 v_0\cdot v_1 + B_1 v_1\cdot v_1) \\
	&+ \int_\Omega (B_0-B_1)v_1\cdot v_1\\
	= &\quad \Re\int _{\partial\Omega} (u_1-u_0)h+\int_D (B_0-B_1)v_1\cdot v_1.
	\end{aligned}
	\end{equation}
	Thus,
	\begin{equation}\label{id1}
	\begin{aligned}
	\Re \delta W&= \int_\Omega B_0 (v_0-v_1)\cdot(v_0-v_1) + \int_D (B_1-B_0)v_1\cdot v_1. \\
	\end{aligned}
	\end{equation}
	Swapping the indices 0 and 1, we obtain a similar identity
	\begin{equation}\label{id2}
	\begin{aligned}
	\Re \delta W=- \int_\Omega B_1 (v_0-v_1)\cdot(v_0-v_1) + \int_D (B_1-B_0)v_0\cdot v_0. 
	\end{aligned}
	\end{equation}
	Here we have used the fact that $\supp\, (B_1-B_0)\subset\bar D$.

	\textbf{Case 1:} \eqref{a0}(i) is satisfied. 
	
	The positivity of $B_1$ and \eqref{id2}  give the second inequality of \eqref{eineq1}.
	Using \eqref{id1} and the triangle inequality, we see that the first inequality of  \eqref{eineq1} follows if $B_1 -B_0$ is positive-definite.  We have $B_1 -B_0 = P+Q$, where
	\[
	P= \begin{pmatrix}M &M\eps_1\\ \eps_1 M& N+ \eps_1 M\eps_1\end{pmatrix}, \quad Q= \begin{pmatrix}0 &\sigma_0^{-1}(\eps_1-\eps_0)\\(\eps_1-\eps_0)\sigma_0^{-1} &  (\eps_1-\eps_0)\sigma_0^{-1}  \eps_1+ \eps_0 \sigma_0^{-1} (\eps_1-\eps_0)\end{pmatrix},
	\]
	and $M=(\sigma_1+\zeta_1)^{-1}-\sigma_0^{-1}$, $N =(\sigma_1-\zeta_1)-\sigma_0$. From  \eqref{a0}(i), it follows that $M$ and $N$ are positive-definite.
	
	For $X=(p,q)^T\in \R^{2n}$, we have 
	\[
	PX\cdot X=M (p+\eps_1 q)\cdot (p+\eps_1 q)
	+ Nq\cdot q
	\]
	and 
	\[
	|QX\cdot X| \le C\|\eps_1-\eps_0\|_{L^\infty(D)} |X|.
	\]
	The required positivity then follows from condition \eqref{delta} for $\delta$ small depending on $\alpha  ,\varrho$.
	
	\textbf{Case 2:} \eqref{a0}(ii) is satisfied.
		
	Similar argument as in the previous case shows that $B_0-B_1$ is positive-definite. Thus, the first inequality of \eqref{eineq2} follows from \eqref{id2}. To prove the second inequality of \eqref{eineq2}, we first deduce from \eqref{id1} and \eqref{id2} that
	\begin{align}\label{new1}
	\int_D(B_0 -B_1) v_1\cdot v_1 &= -\Re \delta W +  \int_\Omega B_0(v_0 - v_1)\cdot (v_0 - v_1) \notag \\
	&\le -\Re \delta W +  C \int_\Omega B_1(v_0 - v_1)\cdot (v_0 - v_1) \notag  \\
	&= -(C+1)\Re \delta W +  C \int_D (B_1-B_0)v_0\cdot v_0.
	\end{align}
	Using \eqref{basic*} and Cauchy-Schwarz we have
	\begin{align}\label{new2}
	-\Re \delta W \le  2(C+1)\int_D (B_0-B_1)v_0\cdot v_0 + \frac{1}{2(C+1)} \int_D (B_0-B_1)v_1\cdot v_1.
	\end{align}
	Eestimating the last term of  \eqref{new2} by  \eqref{new1}, we obtain
	\begin{align}
	-\Re \delta W \le  4(C+1)\int_D (B_0-B_1)v_0\cdot v_0.
	\end{align}
	The proof is complete.
\eproof

Another tool we need is the following Lipschitz propagation of smallness which is a consequence of Theorem~\ref{main-thm-1}.
\begin{pr}\label{small}
Assume that $\partial\Omega\in C^{3,1}$. Let $u_0$ be the solution of \eqref{unpert}. Then for any $a>0$ and every
$x\in\Omega_{4a}$, we have
\begin{equation}
\int_{B_{a}(x)}|\nabla u_0|^{2}\ge C_{a}\int_{\Omega}|\nabla u_0|^{2},\label{lp}
\end{equation}
where $C_{a}$ depends on $\lambda_0$, $M_0$, $\rho_1, M_1$, $\rho_0, K_0$, $|\Omega|$, $a$, and $\|g\|_{L^2(\partial\Omega)}/\|g\|_{H^{-1/2}(\partial\Omega)}$. 
\end{pr}

Before proving this proposition, we first establish a technical lemma.
\begin{lemma}
Let the assumptions described in Proposition~\ref{small} hold. Then for $a>0$ sufficiently small, we have
\begin{equation}\label{claim}
\int_{\Omega\setminus{\Omega_{a/4}}}|\nabla u_0|^2\leq C a^{1/n}\|g\|^2_{L^2(\partial\Omega)}.
\end{equation}
\end{lemma}
\pf. Let $c_0$ be a positive number such that $4c_0<\mbox{\rm dist}(\partial\Omega,\Sigma)$. Denote $\mathcal{U}_{c_0}=\Omega\setminus\Omega_{c_0}$, that is, $\mathcal{U}_{c_0}=\{x\in\Omega\,:\, \mbox{\rm dist}(x,\partial\Omega)\leq c_0\}$. Let $\eta_0$ be a cut-off function, such that $\eta_0\in H^1(\Omega)$ and
\[\eta_0=\left\{\begin{array}{rl}1&\mbox{in }\mathcal{U}_{c_0}\\
0&\mbox{in }\Omega\setminus\mathcal{U}_{2c_0}=\overline{\Omega}_{2c_0}\end{array}\right.\]
Now let $u=u_0\eta_0$, since $\nabla\cdot I_0(\nabla u_0)=0$, we have
\begin{equation*}
\nabla\cdot I_0(\nabla u)= \left[2I_0(\nabla\eta_0)\cdot\nabla u_0+ u_0 \nabla\cdot(I_0(\nabla\eta_0)\right]\chi_{\mathcal{U}_{2c_0}}
\end{equation*}
Moreover, since $u\equiv 0$ in $\overline{\Omega}_{2c_0}$, we have actually that
$\nabla\cdot I_0(\nabla u)=\nabla\cdot A^0_+(\nabla u)$, hence
\[\nabla\cdot (A^0_+\nabla u)=\left[2I_0(\nabla\eta_0)\cdot\nabla u_0+ u_0 \nabla\cdot(I_0(\nabla\eta_0)\right]\chi_{\mathcal{U}_{2c_0}}\]

We now apply Theorem 2.2, Ch. 4, Sect. 4.2.7,  in \cite{n} and get that
\begin{equation}\label{nec}
\|u\|_{H^2(\Omega)}\leq C\left(\|f\|_{L^2(\Omega)}+\|g\|_{H^{1/2}(\partial\Omega)}\right)
\end{equation}
where
\[f=\left[2I_0(\nabla\eta_0)\cdot\nabla u_0+ u_0 \nabla\cdot(I_0(\nabla\eta_0)\right]\chi_{\mathcal{U}_{2c_0}}\]
and $C$ depends on $\lambda_0, M_0, \rho_1, M_1$.
We can estimate
\begin{equation}\label{f}
\begin{aligned}
\|f\|_{L^2(\Omega)}&=\|f\|_{L^2(\mathcal{U}_{2c_0})}\leq C\left(\|\nabla\eta\cdot\nabla u_0\|_{L^2(\mathcal{U}_{2c_0}}+\|u_0\nabla\cdot(I_0(\nabla\eta_0)\|_{L^2(\mathcal{U}_{2c_0}}\right)\\
&\leq C \|u_0\|_{H^1(\Omega)}\leq C\|g\|_{H^{-1/2}(\partial\Omega)}.
\end{aligned}
\end{equation}
Here $C$ depends on $\lambda_0, M_0, \rho_1, M_1,$ and $c_0$.
Combining \eqref{nec}, \eqref{f}, and $\eta_0=1$ in $\mathcal{U}_{c_0}$, we obtain
\begin{equation*}
\|u_0\|_{H^2(\mathcal{U}_{c_0})}\leq C \|g\|_{H^{1/2}(\partial\Omega)}.
\end{equation*}
This estimate, together with the standard energy estimate
\[\|u_0\|_{H^1(\mathcal{U}_{c_0})}\leq C \|g\|_{H^{-1/2}(\partial\Omega)}\]
gives, by interpolation,
\begin{equation}\label{32norm}
\|u_0\|_{H^{3/2}(\mathcal{U}_{c_0})}\leq C \|g\|_{L^2(\partial\Omega)}.
\end{equation}

H\"{o}lder inequality implies
\begin{equation}\label{ars2.21}
\|\nabla u_0\|^2_{L^2\left(\Omega\setminus{\Omega_{a/4}}\right)}\leq \left|\Omega\setminus{\Omega_{a/4}}\right|^{1/n}\|\nabla u_0\|^2_{L^{2n/(n-1)}(\Omega\setminus{\Omega_{a/4}})}.
\end{equation}
By Sobolev inequality and  if $a/4\leq c_0$, we have
\begin{equation}\label{ars2.22}
 \|\nabla u_0\|^2_{L^{2n/(n-1)}(\Omega\setminus{\Omega_{a/4}})}\leq C\|\nabla u_0\|^2_{H^{1/2}(\Omega\setminus{\Omega_{a/4}})}\leq C \|u_0\|_{H^{3/2}(\mathcal{U}_{c_0})}.
\end{equation}
From \eqref{ars2.21}, \eqref{ars2.22} and \eqref{32norm}, and in view of the estimate
\[|\Omega\setminus{\Omega_{a/4}}|\leq C a\]
(see \cite[(A.3)]{ar}), we finally prove \eqref{claim}.\eproof

We now prove Proposition~\ref{small}. 

\smallskip
\noindent{\bf Proof}. It suffices to consider $a$ sufficiently small. Assume $a<h_0$, where $h_0$ is the constant given in  Theorem~\ref{main-thm-1}. To apply Theorem~\ref{main-thm-1}, we choose $h=a/3$ and obtain
\begin{equation}\label{1est}
\|u_0\|_{L^2(\Omega_{a/3})}\le C\|u_0\|^\delta_{L^2(B_a(x))}\|u_0\|^{1-\delta}_{L^2(\Omega)}.
\end{equation}
Using Caccioppoli's inequality and a generalized Poincar\'e inequality \cite[(3.8)]{amr}, we deduce from \eqref{1est} 
\begin{equation*}
\|\nabla u_0\|_{L^2(\Omega_{a/4})}\le C\|\nabla u_0\|^\delta_{L^2(B_a(x))}\|\nabla u_0\|^{1-\delta}_{L^2(\Omega)},
\end{equation*}
i.e.,
\begin{equation}\label{2est}
\frac{\|\nabla u_0\|_{L^2(\Omega_{a/4})}}{\|\nabla u_0\|_{L^2(\Omega)}}\le C\left(\frac{\|\nabla u_0\|_{L^2(B_a(x))}}{\|\nabla u_0\|_{L^2(\Omega)}}\right)^\delta.
\end{equation}

With the help of \eqref{claim} and the estimate
\[
\|g\|^2_{H^{-1/2}(\partial\Omega)}\le C\|\nabla u_0\|^2,
\]
we can see that
\begin{equation}\label{half}
\frac{\|\nabla u_0\|^2_{L^2(\Omega_{a/4})}}{\|\nabla u_0\|^2_{L^2(\Omega)}}=1-\frac{\|\nabla u_0\|^2_{L^2(\Omega\setminus\Omega_{a/4})}}{\|\nabla u_0\|^2_{L^2(\Omega)}}\ge 1-Ca^{1/n}\frac{\|g\|^2_{L^2(\partial\Omega)}}{\|g\|^2_{H^{-1/2}(\partial\Omega)}}\ge \frac 12
\end{equation}
by taking $a$ sufficiently small. Estimate \eqref{lp} follows immediately from \eqref{2est} and \eqref{half}. \eproof

Finally, we provide the proof of the size estimate result. 

\medskip
\noindent{\bf Proof of Theorem \ref{size}}.
We follow  the proof of Theorem 2.1 in \cite[page 61]{ars}. By the estimate for elliptic system with discontinuous coefficients \cite[Theorem~1.1]{ln} and the Poincar\'e inequality, we have the interior estimate
$$
\sup_{D}|\nabla u_0|\leq C\|\nabla u_0\|_{L^2(\Omega)}\leq C \left(\Re W'_0\right)^{1/2},
 $$ 
where $C$ depends on $\lambda_0, M_0, \rho_0, K_0, d_0$. We want to point out that the equation $\nabla\cdot I_0(\nabla u_0)=0$ can be transformed to the elliptic system with required conditions considered in \cite[Theorem~1.1]{ln}. Applying the second inequality of \eqref{eineq1} or \eqref{eineq2}, we obtain
$$
\Re\delta W \leq \kappa_2|D|\left(\sup_D |\nabla u_0|\right)^2\leq C_1^{-1} |D|\Re W'_0, 
$$
which gives a lower bound of $|D|$.

Let $\ell=\min\left(d_0/2,d_1/2\right)$ and cover $D_{d_1}$ with squares $\{Q_k\}_{k=1,\ldots,N}$ of side length $\ell$ and disjoint interiors. It is clear that $N\geq \ell^{-2} |D_{d_1}|\ge\frac{1}{2\ell^2}|D|$. Applying Proposition \ref{small} with $a=\ell/2$ we see that $\int_{Q_k}|\nabla u_0|^2\geq C\Re W'_0$, hence

$$
\int_D |\nabla u_0|^2\geq \sum_{k=1}^N\int_{Q_k}|\nabla u_0|^2\geq N C\Re W'_0 \geq C|D| \Re W'_0.
$$
Combining this with  the first inequality of \eqref{eineq1} or \eqref{eineq2} provides an upper bound of $|D|$.
\eproof

\section*{Acknowledgements}
The research of Francini and Vessella was partly funded by:
(a) Research Project 201758MTR2 of the Italian Ministry of Education, 
University and Research (MIUR)
Prin 2017 “Direct and inverse problems for partial differential 
equations: theoretical aspects and applications”;
(b) GNAMPA of the Italian INdAM – National Institute of High Mathematics 
(grant number not available). Wang was partially supported by MOST 108-2115-M-002-002-MY3 \& 109-2115-M-002-001-MY3.


\begin{thebibliography}{99999999}

\bibitem[AMR]{amr}
G. Alessandrini, A. Morassi, and E. Rosset, \emph{The linear constraints in Poincar\'e and Korn type inequalities}, Forum Math., \textbf{20} (2008), 557-569.

\bibitem[ARRV]{arrv}
G. Alessandrini, L. Rondi, E. Rosset, and S. Vessella, \emph{The stability for the Cauchy problem for elliptic equations}, Inverse problems, \textbf{25} (2009), 123004.

\bibitem[AR]{ar} G Alessandrini and E Rosset, \emph{The inverse conductivity problem with one measurement: bounds on the size of the unknown object}, SIAM J. Appl. Math., \textbf{58} (1998), no. 4, 1060-1071. 

\bibitem[ARS]{ars} G Alessandrini, E Rosset, and J K Seo, \emph{Optimal size estimate for the inverse conductivity problem with one measurement}, Proc. AMS, \textbf{128} (1999), 53-64.


\bibitem[BL]{BL} M. Bellassoued, J. Le Rousseau, \emph{Carleman estimates for elliptic operators with complex coefficients. Part II: Transmission problems}, J. Math. Pures Appl., \textbf{115} (2018), 127-186.

\bibitem[BFV]{bfv}
E. Beretta, E. Francini, and S. Vessella, \emph{Size estimates for the EIT problem with one measurement: the complex case}, Rev. Mat. Iberoam., \textbf{30} (2014), 551-580.

\bibitem[CW]{CW}
C. Carstea and J.-N. Wang, \emph{Propagation of smallness for an elliptic PDE with piecewise Lipschitz coefficients}, J. Differential Equations, \textbf{268} (2020), no. 12, 7609-7628.

\bibitem[CNW]{CNW}
C. Carstea, T. Nguyen, and J.-N. Wang, \emph{Uniqueness estimates for the general complex conductivity equation and their applications to inverse problems}, SIAM J. Math. Anal., \textbf{52} (2020), no. 1, 570-580.

\bibitem[CG]{cg94}
A.V. Cherkaev and L.V. Gibiansky, \emph{Variational principles for complex conductivity, viscoelasticity, and similar problems in media with complex moduli}, J. Math. Phys., \textbf{35(1)} (1994), 127-145.

\bibitem[CGT]{cgt}
F. Colombini, C. Grammatico, D. Tataru, \emph{Strong uniqueness for second order elliptic operators with Gevrey coefficients},  Math. Res. Lett., \textbf{13} (2006), 15-27.

\bibitem[DFLVW]{dflvw}
M. Di Cristo, E. Francini, C.-L. Lin, S. Vessella, and J.-N. Wang, \emph{Carleman estimate for second order elliptic equations with Lipschitz leading coefficients and jumps at an interface}, J. Math. PuresAppl. (9), \textbf{108} (2017), 163-206.


\bibitem[FLVW]{flvw}
E. Francini, C.-L. Lin, S. Vessella, and J.-N. Wang, \emph{Three-region inequalities for the second order elliptic equation with discontinuous coefficients and size estimate},  J. Differential Equations, \textbf{261} (2016), no. 10, 5306-5323.

\bibitem[FVW]{FVW20}
E. Francini, S. Vessella, and J.-N. Wang, \emph{Carleman estimate for complex second order elliptic operators with discontinuous Lipschitz coefficients}, arXiv:2001.04071 [math.AP]

\bibitem[FV]{fv18}
E. Francini and S. Vessella, \emph{Carleman estimates for the parabolic transmission problem and H\"older propagation of smallness across an 
interface}, J. Differential Equations, \textbf{265} (2018), 2375-2430.

\bibitem[Ho1]{ho0} L. H\"ormander, \emph{Linear Partial Differential Operators}, Springer-Verlag, New York, 1969.


\bibitem[Ho3]{ho3} L. H\"ormander, \emph{The Analysis of Linear Partial Differential Operators. Vol. III}, Springer-Verlag, New York, 1985. 

\bibitem[Ho4]{ho4}
L. H\"ormander, \emph{A uniqueness theorem for second order hyperbolic differential equations}, Comm. Partial Differential Equations, \textbf{17} (1992), no. 5-6, 699-714. 

\bibitem[LL]{ll}  J. Le Rousseau and N. Lerner, \emph{Carleman estimates for anisotropic elliptic operators with jumps at an interface}, Analysis \& PDE, \textbf{6} (2013), No. 7, 1601-1648.

\bibitem[LR1]{lr1} J. Le Rousseau and L. Robbiano, \emph{Carleman estimate for elliptic operators with coefficients with jumps at an interface in arbitrary dimension and application to the null controllability of linear parabolic equations}, Arch. Rational Mech. Anal., \textbf{195} (2010), 953-990.

\bibitem[LR2]{lr2} J. Le Rousseau and L. Robbiano, \emph{Local and global Carleman estimates for parabolic operators with coefficients with jumps at interfaces}, Inventiones Math., \textbf{183} (2011), 245-336.

\bibitem[LN]{ln}
Y. Li and L. Nirenberg, \emph{Estimates for elliptic systems from composite material}, Comm. Pure. Appl. Math., \textbf{56} (2003), 0892-0925.

\bibitem[MPH]{mph06}
D. Miklav\u ci\u c, N. Pav\u selj, and F.X. Hart, \emph{Electric Properties of Tissues}, Wiley Encyclopedia of Biomedical Engineering 2006. DOI: 10.1002/9780471740360.ebs0403

\bibitem[M]{mi} K. Miller, \emph{Nonunique continuation for uniformly parabolic and elliptic equations in self-adjoint divergence form with H\"older continuous coefficients}, Arch. Rational Mech. Anal., \textbf{54} (1974), 105-117.

\bibitem[N]{n} J. Ne\u cas, \emph{Direct Methods in the Theory of Elliptic Equations}, Springer-Verlag, Berlin Heidelberg, 2012.

\bibitem[P]{pl} A. Pli\'s, \emph{On non-uniqueness in Cauchy problem for an elliptic second order differential equation}, Bull. Acad. Polon. Sci. S\'er. Sci. Math. Astronom. Phys., \textbf{11} (1963), 95-100. 

\end{thebibliography}
\end{document}